\documentclass[11pt]{amsart}
\usepackage{fullpage, amsmath, amsthm,amsfonts,amssymb,stmaryrd, mathrsfs,amscd}
\usepackage{graphicx}
\usepackage{hyperref}

\usepackage[pdftex]{color}

\newtheorem{thm}{Theorem}[subsection]

\newtheorem{lem-constr}[thm]{Lemma-Construction}

\newtheorem{lem-def}[thm]{Lemma-Definition}
\newtheorem{cor}[thm]{Corollary}

\newtheorem{conjecture}[thm]{Conjecture}

\theoremstyle{definition}

\theoremstyle{definition}

\newtheorem{ex}[thm]{Example}

\newtheorem{rmk}[thm]{Remark}
\theoremstyle{definition}
\newtheorem{dfn}[thm]{Definition}

\numberwithin{equation}{section}

\input xy
\xyoption{all}

\newcommand{\quash}[1]{}  %%Anything in \quash is ignored
\newcommand{\nc}{\newcommand}
\nc{\on}{\operatorname}

\DeclareFontEncoding{OT2}{}{} % to enable usage of cyrillic fonts
\newcommand{\textcyr}[1]{%
 {\fontencoding{OT2}\fontfamily{wncyr}\fontseries{m}\fontshape{n}\selectfont #1}}

\newcommand{\frakg}{{\mathfrak g}}

\newcommand{\frakl}{{\mathfrak l}}
\newcommand{\frakm}{{\mathfrak m}}

\newcommand{\frakA}{{\mathfrak A}}
\newcommand{\frakB}{{\mathfrak B}}

\newcommand{\frakS}{{\mathfrak S}}

\newcommand{\bA}{{\mathbb A}}
\newcommand{\bB}{{\mathbb B}}
\newcommand{\bC}{{\mathbb C}}
\newcommand{\bD}{{\mathbb D}}

\newcommand{\bF}{{\mathbb F}}
\newcommand{\bG}{{\mathbb G}}

\newcommand{\bL}{{\mathbb L}}

\newcommand{\bO}{{\mathbb O}}

\newcommand{\bQ}{{\mathbb Q}}
\newcommand{\bR}{{\mathbb R}}

\newcommand{\bZ}{{\mathbb Z}}

\newcommand{\mA}{{\mathcal A}}
\newcommand{\mB}{{\mathcal B}}

\newcommand{\mE}{{\mathcal E}}
\newcommand{\mF}{{\mathcal F}}
\newcommand{\mG}{{\mathcal G}}
\newcommand{\mH}{{\mathcal H}}

\newcommand{\mL}{{\mathcal L}}

\newcommand{\mO}{{\mathcal O}}

\newcommand{\mS}{{\mathcal S}}

\newcommand{\mZ}{{\mathcal Z}}

\newcommand{\scrS}{{\mathscr S}}

\nc{\al}{{\alpha}} 
\nc{\be}{{\beta}} 
\nc{\ga}{{\gamma}}
\nc{\de}{{\delta}}
\nc{\ep}{{\epsilon}}
\nc{\la}{{\lambda}}

\nc{\Ga}{{\Gamma}} 
\nc{\La}{{\Lambda}}

\nc{\ve}{{\varepsilon}} 

\nc{\Ab}{{\mathbf{Ab}}}   % Category of abelian groups
\nc{\Aff}{{\mathbf{Aff}}}   % Category of affine schemes, affine spaces, etc.
\nc{\Alg}{{\mathbf{Alg}}}  % Category of algebras
\nc{\AlgSp}{{\mathbf{AlgSp}}}  % Category of algebraic spaces
\nc{\Ani}{{\mathbf{Ani}}}   % Category of anima
\nc{\CAlg}{{\mathbf{CAlg}}}   % Category of commutative algebras
\nc{\Cat}{{\mathbf{Cat}}}    % Category of categories
\nc{\Coh}{{\mathbf{Coh}}}        % Category of coherent sheaves
\nc{\Grp}{{\mathbf{Grp}}}        % Category of groups
\nc{\Grpd}{{\mathbf{Grpd}}}        % Category of group groupoids
\nc{\infGrpd}{{\infty\on{-Groupoid}}}  % Category of infinite groupoid
\nc{\LCat}{{\mathbf{LCat}}}  % Category of linear categories
\nc{\Mod}{{\mathbf{Mod}}}
\nc{\Mon}{{\mathbf{Mon}}}   % Category of monoids, monoidal categories
\nc{\Perf}{{\mathbf{Perf}}}    % Category of perfect rings, perfectoid rings, perfect complexes
\nc{\Perv}{{\mathbf{Perv}}}  % Category of perverse sheaves
\nc{\Pre}{{\mathbf{Pre}}}     % Category of presentable categories
\nc{\PStk}{{\mathbf{PStk}}}  % Category of prestacks
\nc{\QCoh}{{\mathbf{QCoh}}}  % Category of quasi-coherent sheaves
\nc{\Rep}{{\mathbf{Rep}}}   % Category of representations
\nc{\sCat}{{\mathbf{sCat}}}   % Category of small categories
\nc{\sLCat}{{\mathbf{sLCat}}}  % Category of small linear categories
\nc{\Sch}{{\mathbf{Sch}}}   % Category of schemes  
\nc{\Sets}{\mathbf{Sets}}   % Category of sets
\nc{\Shv}{{\mathbf{Shv}}}       % Category of sheaves
\nc{\Spc}{{\mathbf{Spc}}}  % Category of spaces
\nc{\Stk}{{\mathbf{Stk}}}  % Category of stacks
\nc{\SymMon}{{\mathbf{SymMon}}} % Category of symmetric monoidal categories
\nc{\Top}{{\mathbf{Top}}}   % Category of topological spaces
\nc{\Tor}{{\mathbf{Tor}}}  % torsor
\nc{\Vect}{{\mathbf{Vect}}} % Category of vector spaces
\nc{\fgMod}{\Mod^{\on{f.g.}}}
\nc{\fgRep}{\Rep_{\on{f.g.}}}

\nc{\alg}{{\on{alg}}}
\nc{\Aut}{{\on{Aut}}}
\nc{\aut}{{\on{aut}}}
\nc{\Bcrs}{\on{B}_{\cris}}
\nc{\BdR}{\on{B}_{\dR}}
\nc{\Bst}{\on{B}_{\st}}
\nc{\cha}{{\on{char}}}       % charateristic
\nc{\Cha}{{\on{Char}}}       % charateristic
\nc{\coh}{{\on{coh}}}
\nc{\coker}{{\on{coker}}}
\nc{\conv}{{\on{conv}}}     % convolution
\nc{\Corr}{{\on{Corr}}}       % correspondence
\nc{\corr}{{\on{corr}}} 
\nc{\Der}{{\on{Der}}}
\nc{\der}{{\on{der}}}
\nc{\diag}{{\on{diag}}}
\nc{\End}{{\on{End}}}
\nc{\et}{{\on{et}}}
\nc{\ext}{{\on{ext}}}
\nc{\Ext}{{\on{Ext}}}
\nc{\Fil}{{\on{Fil}}}
\nc{\Fr}{{\on{Frob}}}
\nc{\Fun}{{\on{Fun}}}  % function, functor
\nc{\Gal}{{\on{Gal}}}  % Galois
\nc{\gr}{{\on{gr}}}    % graded, grassmannian
\nc{\Gr}{{\on{Gr}}}   % graded, grassmannian
\nc{\heart}{\ensuremath\heartsuit}
\nc{\Hom}{{\on{Hom}}}
\nc{\IC}{{\on{IC}}}    % Intersection coh sheaf
\nc{\id}{{\on{id}}}
\nc{\Id}{{\on{Id}}}
\nc{\ind}{{\on{ind}}}
\nc{\Ind}{{\on{Ind}}}
\nc{\Int}{{\on{Int}}}
\nc{\inv}{{\on{Inv}}}
\nc{\Iso}{{\on{Isom}}}
\nc{\Lie}{{\on{Lie}}}
\nc{\Ker}{{\on{Ker}}}
\nc{\Map}{{\on{Map}}}
\nc{\mult}{{\on{mult}}}    % Multiplication, multiplicity 
\nc{\Nm}{{\on{Nm}}}      % Norm
\nc{\Ql}{{\overline{\bQ}_\ell}}
\nc{\out}{\on{out}}
\nc{\Out}{\on{Out}}
\nc{\Orb}{{\on{Orb}}}
\nc{\ord}{{\on{ord}}}
\nc{\Pic}{{\on{Pic}}}
\nc{\pr}{{\on{pr}}}         % Projection
\nc{\Pro}{{\on{Pro}}}    % Pro-objects
\nc{\proj}{{\on{proj}}}
\nc{\Proj}{{\on{Proj}}} 
\nc{\pt}{{\on{pt}}}        % point
\nc{\res}{{\on{res}}}
\nc{\Res}{{\on{Res}}}
\nc{\RHom}{\on{RHom}}
\nc{\sgn}{{\on{sgn}}}
\nc{\sha}{{\mbox{\textcyr{Sh}}}}  
\nc{\Spa}{{\on{Spa}}}
\nc{\Spd}{{\on{Spd}}}
\nc{\Spec}{{\on{Spec}}}
\nc{\Spf}{{\on{Spf}}}
\nc{\std}{\on{std}}   % Standard
\nc{\Std}{\on{Std}}  % Standard
\nc{\sym}{\on{sym}}
\nc{\Sym}{\on{Sym}}  % symmetric
\nc{\tr}{{\on{tr}}}   % Trace
\nc{\Tr}{{\on{Tr}}}   % Trace
\nc{\xch}{\mathbb{X}^\bullet}      % Weight lattice
\nc{\xcoch}{\mathbb{X}_\bullet}  % Coweight lattice

\nc{\Bun}{{\on{Bun}}}
\nc{\Fl}{{\mF\ell}}       % Flag variety

\nc{\Loc}{{\on{Loc}}}
\nc{\Hk}{{\on{Hk}}}
\nc{\lochk}{{\Hk^{\loc}}}
\nc{\Sh}{{\on{Sh}}}   % Shimura variety
\nc{\Sht}{{\on{Sht}}}  % Moduli of Shtukas

\nc{\GL}{{\on{GL}}}
\nc{\GSp}{{\on{GSp}}} 
\nc{\GSpin}{{\on{GSpin}}} 
\nc{\GU}{{\on{GU}}} 
\nc{\PGL}{{\on{PGL}}}
\nc{\SL}{{\on{SL}}}
\nc{\Sp}{{\on{Sp}}}
\nc{\Spin}{{\on{Spin}}}  
\nc{\SU}{{\on{SU}}} 
\nc{\SO}{{\on{SO}}}
\nc{\gl}{{\frakg\frakl}}
\nc{\CT}{{}^cT}  
\nc{\CB}{{}^cB}  
\nc{\CG}{{}^cG}     % C-group
\nc{\LG}{{}^LG}       % L-group

\nc{\ab}{{\on{ab}}}    % Abelianization
\nc{\ad}{{\on{ad}}}    % Adjoint          
\nc{\Ad}{{\on{Ad}}}   % Adjoint   
\nc{\adm}{{\on{adm}}}   % Admissible
\nc{\aff}{{\on{aff}}}
\nc{\an}{{\on{an}}}   % analytification
\nc{\Ch}{{\on{Ch}}}
\nc{\CH}{{\on{CH}}}
\nc{\cl}{{\on{cl}}}
\nc{\cons}{{\on{cons}}}  % constructible
\nc{\cris}{{\on{cris}}}
\nc{\crys}{{\on{crys}}}
\nc{\cycl}{{\on{cycl}}}
\nc{\DCA}{\mathbf{DCA}}
\nc{\DerAff}{{\on{DerAff}}}
\nc{\disc}{{\on{disc}}}
\nc{\dR}{{\on{dR}}}
\nc{\Iw}{{\on{Iw}}}  % Iwahori
\nc{\loc}{{\on{loc}}}
\nc{\op}{{\on{op}}}
\nc{\pf}{{\on{pf}}}
\nc{\reg}{{\on{reg}}}
\nc{\Tate}{{\on{Tate}}}
\nc{\s}{{\on{sc}}}        % simply connected
\nc{\Sat}{{\on{Sat}}}   % Satake
\nc{\sep}{{\on{sep}}}
\nc{\Sph}{{\on{Sph}}}  % Spherical
\nc{\st}{{\on{st}}}    %  Steinberg, stable
\nc{\St}{{\on{St}}}   %  Steinberg, stable
\nc{\ta}{{\on{tame}}}
\nc{\un}{{\on{unip}}}
\nc{\ur}{{\on{ur}}}
\nc{\wt}{{\on{wt}}}

\nc{\rec}{{\on{rec}}}
\nc{\vect}{{\on{vec}}}
\nc{\kot}{{\frakB}}

\newcommand{\xdashrightarrow}[2][]{\ext@arrow 3359 \rightarrowfill@@{#1}{#2}}
\def\rightarrowfill@@{\arrowfill@@\relax\relbar\shortrightarrow}
\def\arrowfill@@#1#2#3#4{%
  $\m@th\thickmuskip0mu\medmuskip\thickmuskip\thinmuskip\thickmuskip
   \relax#4#1
   \xleaders\hbox{$#4#2$}\hfill
   #3$%
}

\begin{document}

\title{Arithmetic and Geometric Langlands Program}
\author{Xinwen Zhu}

\maketitle

\begin{abstract}
We explain how the geometric Langlands program inspires some recent new prospectives of classical arithmetic Langlands program and leads to the solutions of some problems in arithmetic geometry.
\end{abstract}

\tableofcontents

The classical Langlands program, originated by Langlands in 1960s \cite{Langlands-to-Weil}, systematically studies reciprocity laws in the framework of representation theory.
Very roughly speaking, it predicts the following deep relations between number theory and representation theory.
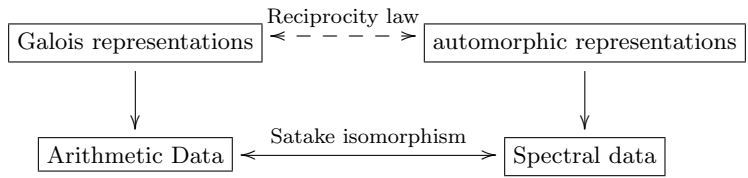
\[
\xymatrix{
 \framebox{\txt{\footnotesize Galois representations}}  \ar@{<-->}^-{\txt{\scriptsize Reciprocity law}}[rr]  \ar[d]&&  \framebox{\txt{\footnotesize automorphic representations}}\ar[d] \\
\framebox{\footnotesize Arithmetic Data} \ar@{<->}^-{\txt{\scriptsize Satake isomorphism}}[rr] &&  \framebox{\footnotesize Spectral data} 
}
\]
A special case of this correspondence, known as the Shimura-Tanniyama-Weil conjecture, implies Fermat's last theorem (see \cite{Wiles-Fermat}).

The geometric Langlands program \cite{Laumon-Geometric-Langlands}, initiated by Drinfeld and Laumon, arose as a generalization of  Drinfeld's approach (\cite{Drinfeld-GL2}) to the global Langlands correspondence for $\GL_2$ over function fields. In the geometric theory, the fundamental object to study shifts from the space of automorphic forms of a reductive group $G$ to the category of sheaves on the moduli space of $G$-bundles on an algebraic curve.

For a long time, developments of the geometric Langlands were inspired by problems and techniques from the classical Langlands, with another important source of inspiration from quantum physics. The basic philosophy is known as categorification/geometrization.
In recent years, however, the geometric theory has found fruitful applications to the classical Langlands program and some related arithmetic problems. Traditionally, one applies Grothendieck's sheaf-to-function dictionary to
``decategorify" sheaves studied in geometric theory to obtain functions studied in arithmetic theory. This was used in Drinfeld's approach to the Langlands correspondence for $\GL_2$, as mentioned above. Another celebrated example is Ng\^o's proof of the fundamental lemma (\cite{Ngo-Fundamental-Lemma}). 
In recent years, there appears another passage from the geometric theory to the arithmetic theory, again via a trace construction,  but is of different nature and is closely related to ideas from physics. 
V. Lafforgue's work on the global Langlands correspondence over function fields (\cite{Lafforgue-Automorphic-to-Galois}) essentially (but implicitly) used this idea.

In this survey article, we review (a small fraction of) the developments of the geometric Langlands, and discuss some recent new prospectives of the classical Langlands inspired by the geometric theory, which in turn lead solutions of some concrete arithmetic problems.
The following diagram can be regarded as a road map. 
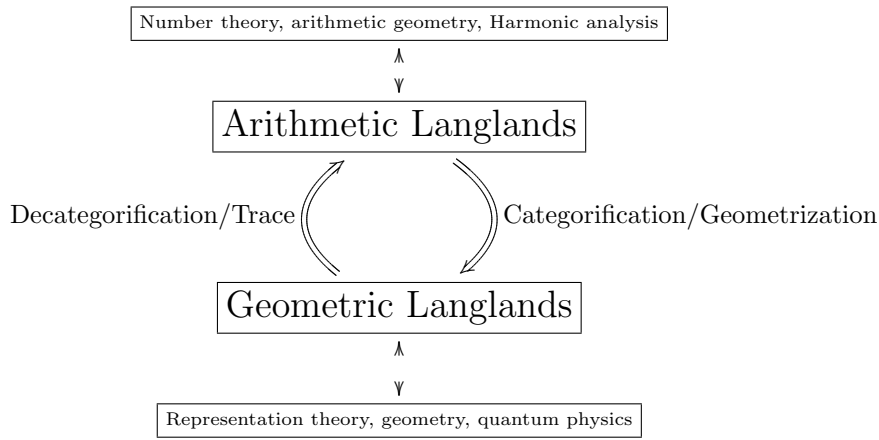
\[
\xymatrix@R=1.5pc{
\framebox{\tiny Number theory, arithmetic geometry,  Harmonic analysis} \ar@{<-->}[d]\\
\framebox{\Large Arithmetic Langlands}\ar^-{\text{\small Categorification/Geometrization}}@{=>}@/^3pc/[dd]\\
\\
\ar^-{\text{\small Decategorification/Trace}}@{=>}@/^3pc/[uu]\framebox{\Large Geometric Langlands} \\ 
\framebox{\tiny Representation theory, geometry, quantum physics} \ar@{<-->}[u]
}
\]

\subsection*{Notations and conventions}
We use the following notations throughout this article.
For a field $F$,  let $\Ga_{\widetilde F/F}$ be the Galois group of a Galois extension $\widetilde F/F$. Write $\Ga_F=\Ga_{\overline F/F}$, where $\overline F$ is a separable closure of $F$. Often in the article $F$ will be either a local or a global field. In this case, let $W_F$ denote the Weil group of $F$. Let $\cycl$ denote the cyclotomic character.

For a group $A$ of multiplicative type over a field $F$, let $\xch(A)=\Hom(A_{\overline F},\bG_m)$ denote the group of its characters, and $\xcoch(A)=\Hom(\bG_m, A_{\overline {F}})$ the group of its cocharacters.

For a prime $\ell$, let $\La$ be $\bF_\ell$, $\bZ_\ell$, $\bQ_\ell$ or a finite (flat) extension of such rings. It will serve as the coefficient ring of our sheaf theory.

\subsection*{Acknowledgement} The author would like to thank all of his collaborators,  without whom many works reported in this survey article would not be possible. Besides, he would also like to thank Edward Frenkel, Dennis Gaitsgory, Xuhua He, Michael Rapoport, Peter Scholze for teaching him and for discussions on various parts of the Langlands program over years. The work is partially supported by NSF under agreement Nos. DMS-1902239 and a Simons fellowship.

\section{From classical to geometric Langlands correspondence}
In this section, we review some developments of the geometric Langlands theory inspired from the classical theory, with another important source of inspiration from quantum physics. The basic idea is categorification/geometrization, which is a process of replacing set-theoretic statements with categorical analogues
\begin{equation}\label{categorification}
\text{Numbers}\dashrightarrow \text{Vector spaces} \dashrightarrow \text{Categories} \dashrightarrow \text{2-Categories} \dashrightarrow \cdots. 
\end{equation}
We illustrate this process by some important examples.

\subsection{The geometric Satake}\label{S: geometric Satake}
The starting point of the Langlands program is  (Langlands' interpretation of) the Satake isomorphism, in which the Langlands dual group appears mysteriously.
Similarly, the starting point of the geometric Langlands theory is the geometric Satake equivalence, which is a tensor equivalence between the category of perverse sheaves on the (spherical) local Hecke stack of a connected reductive group and the category of finite dimensional algebraic representations of its dual group. This is a vast generalization of the classical Satake isomorphism (via the sheaf-to-function dictionary), and arguably gives a conceptual explanation why the Langlands dual group (in fact the $C$-group) should appear in the Langlands correspondence.

We follow \cite[Sect. 1.1]{Zhu-Integral-Satake} for notations and conventions regarding dual groups. For a connected reductive group $G$ over a field $F$, let $(\hat G, \hat B, \hat T ,\hat{e})$ be a pinned Langlands dual group of $G$ over $\bZ$. There is a finite Galois extension $\widetilde F/F$, and a natural injective map $\xi:\Ga_{\widetilde F/F}\subset  \Aut(\hat G, \hat B, \hat T, \hat e)$,  induced by the action of $\Ga_F$ on the root datum of $G$. Let $\LG=\hat{G}\rtimes\Ga_{\widetilde F/F}$ denote the usual $L$-group of $G$, and $\CG=\hat{G}\rtimes (\bG_m\times\Ga_{\widetilde F/F})$ the group defined in \cite{Zhu-Integral-Satake}, which is isomorphic to the $C$-group of $G$ introduced by Buzzard-Gee. We write $d: \CG\to \bG_m\times\Ga_{\widetilde F/F}$ for the projection with kernel $\hat{G}$.

Let $F$ be a non-archimidean local field with $\mO$ its ring of integers and $k=\bF_q$ its residue field. I.e. $F$ is a finite extension of $\bQ_p$ or is isomorphic to $\bF_q(\!(\varpi)\!)$.
Let $\sigma$ be the \emph{geometric} $q$-Frobenius of $k$.
Assume that $G$ can be extended to a connected reductive group over $\mO$ (such $G$ is called unramified) and we fix such an extension so we have $G(\mO)\subset G(F)$, usually called a hyperspecial subgroup of $G(F)$. With a basis of open neighborhoods of the unit given by finite index subgroups of $G(\mO)$,
the group $G(F)$ is a locally compact topological group. 
The classical spherical Hecke algebra is the space of compactly supported $G(\mO)$-bi-invariant $\bC$-valued functions on $G(F)$, equipped with the convolution product
\begin{equation}
\label{E:conv prod alge}
(f*g)(x)=\int_{G(F)}f(y)g(y^{-1}x)dy,
\end{equation}
where $dy$ is the Haar measure on $G(F)$ such that $G(\mO)$ has volume $1$.
Note that if both $f$ and $g$ are $\bZ$-valued, so is $f*g$. Therefore, the subset $H^{\cl}_{G(\mO)}$ of $\bZ$-valued functions form a $\bZ$-algebra\footnote{Here $(-)^{\cl}$ stands for the classical Hecke algebra, as opposed to the derived Hecke algebra mentioned in \eqref{E: Hecke alg isom}.}.

On the dual side,  under the unramifiedness assumption, $\Ga_{\widetilde F/F}$ is a finite cyclic group generated by $\sigma$. Note that $\hat{G}$ acts on $\CG|_{d=(q,\sigma)}$, the fiber of $d$ at $(q,\sigma)\in \bG_m\times \Ga_{\widetilde F/F}$, by conjugation.
Then the classical Satake isomorphism establishes a canonical isomorphism of $\bZ[q^{-1}]$-algebras
\begin{equation}\label{Sat isom}
\Sat^{\cl}: \bZ[q^{-1}][\CG|_{d=(q,\sigma)}]^{\hat{G}}\cong H^{\cl}_{G(\mO)}\otimes \bZ[q^{-1}].
\end{equation}

\begin{rmk}
In fact, as explained in \cite{Zhu-Integral-Satake}, there is a Satake isomorphism over $\bZ$ (without inverting $q$), in which the $C$-group $\CG$ is replaced by certain affine monoid containing it as the group of invertible elements.
On the other hand, if we extend the base ring from $\bZ[q^{-1}]$ to $\bZ[q^{\pm \frac{1}{2}}]$, one can rewrite \eqref{Sat isom} as an isomorphism 
\begin{equation}\label{Sat isom classical}
 \bZ[q^{\pm \frac{1}{2}}][\hat{G}\sigma]^{\hat{G}}\cong H^{\cl}_{G(\mO)}\otimes \bZ[q^{\pm \frac{1}{2}}],
\end{equation}
where $\hat{G}$ acts on $\hat{G}\sigma\subset {}^LG$ by the usual conjugation  (e.g. see \cite{Zhu-Integral-Satake} for the discussion).  This is the more traditional formulation of the Satake isomorphism, which is slightly non-canonical, but suffices for many applications.
\end{rmk}

In the geometric theory, where instead of thinking $G(F)$ as a topological group and considering the space of $G(\mO)$-bi-invariant compactly supported functions on it, one regards $G(F)$ as certain algebro-geometric object and studies the category of $G(\mO)$-bi-equivariant sheaves on it. 
In the rest of the section, we allow $F$ to be slightly more general. Namely, we assume that
$F$ is a local field complete with respect to a discrete valuation, with ring of integers $\mO$ and a \emph{perfect} residue field $k$ of characteristic $p>0$ \footnote{If $\cha F=\cha k$ (the equal characteristic case), this assumption on $k$ is not necessary. We impose it here to have a uniform treatment of both equal and mixed characteristic (i.e. $\cha F\neq \cha k$) cases. For the same reason, we work with perfect algebraic geometry below even in equal characteristic.}. Let $\varpi\in\mO$ be a uniformizer. 

We work in the realm of perfect algebraic geometry. 
Recall that a $k$-algebra $R$ is called perfect if the Frobenius endomorphism $R\to R, \ r\mapsto r^p$ is a bijection. Let $\Aff_k^\pf$ denote the category of perfect $k$-algebras. By a perfect presheaf (or more generally a perfect prestack), we mean a functor from $\Aff_k^\pf$ to the category $\Sets$ of sets (or more generally a functor from $\Aff_k^\pf$ to the $\infty$-category $\Spc$ of spaces).  
Many constructions in usual algebraic geometry work in this setting. E.g. one can endow $\Aff_k^\pf$ with Zariski, \'etale or fpqc topologies as usual and talk about corresponding sheaves and stacks. One can then define perfect schemes, perfect algebraic spaces, perfect algebraic stacks, etc., as sheaves (stacks) with certain properties. It turns out that the category of perfect schemes/algebraic spaces defined this way is equivalent to the category of perfect schemes/algebraic spaces in the usual sense.
Some foundations of perfect algebraic geometry can be found in \cite[Appendix A]{Zhu-Satake-Mixed}, \cite{Bhatt-Scholze-Witt} and \cite[\S A.1]{Xiao-Zhu-Cycle}.

For a perfect $k$-algebra $R$, let $W_\mO(R)$ denote the ring of Witt vectors in $R$ with coefficient in $\mO$. If $\cha F=\cha k$, then $W_\mO(R)\simeq R[[\varpi]]$. If $\cha F\neq \cha k$, see \cite[Sect. 0.5]{Zhu-Satake-Mixed}.
If $R=\overline k$, we write $W_\mO(\overline k)$ by $\mO_{\breve F}$ and $W_{\mO}(\overline k)[1/\varpi]$ by $\breve F$.
We write $D_R=\Spec W_\mO(R)$ and $D_R^*=\Spec W_\mO(R)[1/\varpi]$
which are thought as a family of (punctured) discs parameterized by $\Spec R$. 

We denote by $L^+G$ (resp. $LG$) the \emph{jet group} (resp. \emph{loop group}) of $G$. As presheaves on $\Aff_k^\pf$,
\begin{equation*}\label{E:jet and loop}
L^+G(R)=G(W_\mO(R)),\quad LG(R)=G(W_\mO(R)[1/\varpi]).
\end{equation*}
Note that $L^+G(k)=G(\mO)$ and $LG(k)=G(F)$. Let
\[
\Hk_G:=L^+G\backslash LG/L^+G
\]
be the \'etale stack quotient of $LG$ by the left and right $L^+G$-action, sometimes called the (spherical) local Hecke stack of $G$. As a perfect prestack, it sends $R$ to triples $(\mE_1,\mE_2,\beta)$, where $\mE_1,\mE_2$ are two $G$-torsors on $D_R$, and $\beta: \mE_1|_{D_R^*}\simeq \mE_2|_{D_R^*}$ is an isomorphism. 

For $\ell\neq p$,
the modern developments of higher category theory allow one to define the $\infty$-category of \'etale $\bF_\ell$-sheaves on any prestack (e.g. see \cite{Hemo-Zhu-Unipotent}). 
In particular, for $\La=\bF_\ell, \bZ_\ell, \bQ_\ell$ (or finite extension of these rings), it is possible to define the $\infty$-category $\Shv(\Hk_G,\La)$ of $\La$-sheaves on $\Hk_{G}$, which is the categorical analogue of the space of $G(\mO)$-bi-invariant functions on $G(F)$. But without knowing some geometric properties of $\Hk_{G}$, very little can be said about $\Shv(\Hk_G,\La)$. The crucial geometric input is the following theorem.
\begin{thm}\label{affGrass}
Let $\Gr_G:=LG/L^+G$ be the \'etale quotient of $LG$ by the (right) $L^+G$-action, which admits the left $L^+G$-action. Then $\Gr_G$ can be written as an inductive limit of $L^+G$-stable subfunctors $\underrightarrow\lim X_i$, with each $X_i$ being a perfect projective variety and $X_i\to X_{i+1}$ being a closed embedding. 
\end{thm}
The space $\Gr_G$ is usually called the affine Grassmannian of $G$. See  \cite{Beauville-Laszlo, Faltings-Loop} for the equal characteristic case and \cite{Zhu-Satake-Mixed, Bhatt-Scholze-Witt} for the mixed characteristic case, and see \cite{Zhu-PCMI, Zhu-CDM} for examples of closed subvarieties in $\Gr_G$.  The theorem allows one to define the category of constructible and perverse sheaves on $\Hk_G$, and to formulate the geometric Satake, as we discuss now.

First, the (left) quotient by $L^+G$-action induces a map $\Gr_G\to\Hk_G$. Roughly speaking, a sheaf on $\Hk_G$ is perverse (resp. constructible) if its pullback to $\Gr_G$ comes from a perverse (resp. constructible) sheaf on some $X_i$. Then insdie $\Shv(\Hk_G,\La)$ we have the categories $\Perv(\Hk_G,\La)\subset\Shv_{c}(\Hk_G,\La)$ of perverse and constructible sheaves on $\Hk_G$. 
They can be regarded as categorical analogues of the space of $G(\mO)$-bi-invariant compactly supported functions on $G(F)$.  In addition, $\Perv(\Hk_G,\La)$ is an abelian category, semisimple if $\La$ is a field of characteristic zero, called the Satake category of $G$. For simplicity, we assume that $\La$ is a field in the sequel.\footnote{The formulation for $\La=\bZ_\ell$ is slightly more complicated, as the right hand side of \eqref{E: conv product of sheaves} may not be perverse when $\mA$ and $\mB$ are perverse.}

There is also a categorical analogue of the convolution product \eqref{E:conv prod alge}. Namely, there is the convolution diagram
\[
 \Hk_G\times\Hk_G\xleftarrow{\pr}  L^+G\backslash LG\times^{L^+G} LG/L^+G\xrightarrow{m} \Hk_G,
\]
and the convolution of two sheaves $\mA,\mB\in\Shv(\Hk_G,\La)$ is defined as
\begin{equation}\label{E: conv product of sheaves}
\mA\star\mB:= m_!\pr^*(\mA\boxtimes \mB).
\end{equation} 
This convolution product makes $\Shv(\Hk_G,\La)$ into a monoidal $\infty$-category containing $\Perv(\Hk_G,\La)\subset \Shv_{c}(\Hk_G,\La)$ as monoidal subcategories. 

\begin{rmk}
The above construction of the Satake category as a monoidal category is essentially equivalent to the more traditional approach, in which the Satake category is defined as the category of $L^+G$-equivariant perverse sheaves on $\Gr_G$ (e.g. see \cite{Zhu-PCMI} for an exposition). 
\end{rmk}

Let $\Coh(\bB\hat{G}_\La)^{\heart}$ denote the abelian monoidal category of coherent sheaves on the classifying stack $\bB\hat{G}_\La$ over $\La$\footnote{In the dual group side, we always work in the realm of usual algebraic geometry, so $\bB\hat{G}$ is an Artin stack in the usual sense.}, which is equivalent to the category of
algebraic representations of $\hat{G}$ on finite dimensional $\La$-vector spaces. This following theorem is usually known as the geometric Satake equivalence.
\begin{thm}\label{intro:geomSat}
There is a canonical equivalence of monoidal abelian categories
\[
\Sat_G:  \Coh(\bB\hat{G}_\La)^{\heart}\cong \Perv(\Hk_G\otimes \overline k,\La).
\]
\end{thm}
Geometric satake is really one of the cornerstones of the geometric Langlands program, and has been found numerous applications to representation theory, mathematical physics, and (arithmetic) algebraic geometry.
When $F=k(\!(\varpi)\!)$, this theorem grew out of works of Lusztig, Ginzburg, Beilinson-Drinfeld and Mirkovi\'{c}-Vilonen (cf. \cite{Lusztig-q-analog, Beilinson-Drinfeld-Hitchin, Mirkovic-Vilonen-GeoSat}). In mixed characteristic, it was proved in \cite{Zhu-Satake-Mixed, Yu-Integral-Satake}, with the equal characteristic case as an input, and in \cite{Fargues-Scholze-Geometrization} by mimicking the strategy in equal characteristic. 
We conclude this subsection with a few remarks.

\begin{rmk}
\label{R: birth of the dual group}
(1) As mentioned before, the geometric Satake can be regarded as the conceptual definition of the Langlands dual group $\hat G$ of $G$, namely as the Tannakian group of the Tannakian category $\Perv(\Hk_G\otimes \overline k,\La)$. In addition, as explained in \cite{Yun-Zhu-Integral-Homology, Zhu-Ramified-Sat}, the group $\hat G$ is canonically equipped with a pinning $(\hat B, \hat T, \hat e)$. In the rest of the article, by the pinned Langlands dual group $(\hat{G},\hat{B},\hat{T},\hat{e})$ of $G$, we mean the quadruple defined by the geometric Satake.

(2) For arithmetic applications, one needs to understand the $\Ga_k$-action on $\Perv(\Hk_G\otimes\overline k,\La)$ in terms of the dual group side.  It turns out that such action is induced by an action of $\Ga_k$ on $\hat{G}$, preserving $(\hat{B},\hat{T})$ but not $\hat{e}$. See \cite{Zhu-Ramified-Sat, Zhu-PCMI}, or \cite{Zhu-CDM} from the motivic point of view. This leads the appearance of the group $\CG$. See \cite{Zhu-Ramified-Sat, Zhu-PCMI, Zhu-Integral-Satake} for detailed discussions.

(3) There is also  the derived Satake equivalence \cite{Bezrukavnikov-Finkelberg-Derived-Satake}, describing $\Shv_{c}(\Hk_G\otimes\overline k,\La)$ in terms of the dual group, at least when $\La$ is a field of characteristic zero. We mention that the category in the dual side is not the derived category of coherent sheaves on $\bB\hat{G}_\La$.

(4) In fact, for many applications, it is important to have a family version of the geometric Satake. For a (non-empty) finite set $S$,
there is a local Hecke stack $\Hk_{G,D^S}$ over $D^S$, the self-product of the disc $D=\Spec \mO$ over $S$, which roughly speaking classifies quadruples $(\{x_s\}_{s\in S}, \mE,\mE',\beta)$, where $\{x_s\}_{s\in S}$ is an $S$-tuple of points of $D$, $\mE$ and $\mE'$ are two $G$-torsors on $D$, and $\beta$ is an isomorphism between $\mE$ and $\mE'$ on $D-\cup_{s} \{x_s\}$. In equal characteristic, one can regard $D$ as the formal disc at a $k$-point of an algebraic curve $X$ over $k$ and $\Hk_{G,D^S}$ is the restriction along $D^S\to X^S$ of the stack 
\[
\Hk_{G, X^S}=(L^+G)_{X^S}\backslash (LG)_{X^S}/(L^+G)_{X^S},
\] 
where $(LG)_{X^S}$ and $(L^+G)_{X^S}$ are family versions of $LG$ and $L^+G$ over $X^S$ (e.g. see \cite[Sect. 3.1]{Zhu-PCMI} for precise definitions). In mixed characteristic, the stack $\Hk_{G,D^S}$ (and in fact $D^S$ itself) does not live in the world of (perfect) algebraic geometry, but rather in the world of perfectoid analytic geometry as developed by Scholze  (see \cite{Scholze-Berkeley, Fargues-Scholze-Geometrization}). In both cases, one can consider certain category $\Perv^{\on{ULA}}(\Hk_{G,D^S}\otimes\overline k,\La)$ of (ULA) perverse sheaves on $\Hk_{G,D^S}\otimes\overline k$. In addition, for a map $S\to S'$ of finite sets, restriction along $\Hk_{G,D^{S'}}\to\Hk_{G,D^S}$ gives a functor $\Perv^{\on{ULA}}(\Hk_{G,D^S}\otimes\overline k,\La)\to \Perv^{\on{ULA}}(\Hk_{G,D^{S'}}\otimes\overline k,\La)$\footnote{Such restriction functor defines the so-called fusion product,  a key concept in the geometric Satake equivalence. The terminology ``fusion" originally comes from conformal field theory.}. We refer to the above mentioned references for details.

On the other hand, let $\hat{G}^S$ be the $S$-power self-product of $\hat{G}$ over $\La$. Then for $S\to S'$, restriction along $\bB \hat{G}^{S'}\to \bB \hat{G}^S$ gives a functor $\Coh(\bB \hat{G}_{\La}^S)^\heart\to \Coh(\bB \hat{G}_{\La}^{S'})^\heart$.
Now a family version of the geometric Satake gives a system of functors
\begin{equation}\label{factorizable Sat}
\Sat_S: \Coh( \bB \hat{G}_{\La}^S)^\heart\to \Perv^{\on{ULA}}(\Hk_{G,D^S}\otimes\overline k,\La),
\end{equation}
compatible with restriction functors on both sides induced by maps between finite sets (see \cite{Gaitsgory-deJong, Zhu-PCMI}).

(5) For applications, it is important to have the geometric Satake in different sheaf theoretic contents over different versions of local Hecke stacks. Besides the above mentioned ones, we also mention a $D$-module version \cite{Beilinson-Drinfeld-Hitchin}, and an arithmetic $D$-module version \cite{Xu-Zhu-Bessel}.
\end{rmk}

\subsection{Tamely ramified local geometric Langlands correspondence}\label{Tamely local}
We first recall the classical theory.  Assume that $F$ is a finite extension of $\bQ_p$ or is isomorphic to $\bF_q(\!(\varpi)\!)$, and for simplicity assume that $G$ extends to a connected reductive group over $\mO$. (In fact, results in the subsection hold in appropriate forms for quasi-split groups that are split over a tamely ramified extension of $F$.) In addition, we fix a pinning $(B,T,e)$ of $G$ over $\mO$.

The classical local Langlands program aims to classify (smooth) irreducible representations of $G(F)$ (over $\bC$) in terms of Galois representations.
From this point of view, the Satake isomorphism \eqref{Sat isom} gives a classification of irreducible unramified representations of $G(F)$, i.e. those that have non-zero vectors fixed by $G(\mO)$, as such representations are in one-to-one correspondence with simple modules of $H^{\cl}_{G(\mO)}\otimes\bC$, which via the Satake isomorphism \eqref{Sat isom} are parameterized by semisimple $\hat{G}$-conjugacy classes in $\CG$. (For an irreducible unramified representation $\pi$, the corresponding $\hat{G}$-conjugacy class in $\CG$ is usually called the Satake parameter of $\pi$.)

The next important class of irreducible representations are those that have non-zero vectors fixed by an Iwahori subgroup $G(F)$. 
For example, under the reduction mod $\varpi$ map $G(\mO)\to G(k)$, the preimage $I$ of $B(k)\subset G(k)$ is an Iwahori subgroup of $G(F)$. As in the unramified case, the $\bZ$-valued $I$-bi-invariant functions form a $\bZ$-algebra $H^{\cl}_{I}$ with multiplication given by convolution \eqref{E:conv prod alge} (with the Haar measure $dg$ chosen so that the volume of $I$ is one), and the set of irreducible representations of $G(F)$ that have non-zero $I$-fixed vectors are in bijection with the set of simple $(H^{\cl}_I\otimes\bC)$-modules.
Just as the Satake isomorphism, Kazhdan-Lusztig gave a description of $H^{\cl}_{I}\otimes\bC$ in terms of geometric objects associated to $\hat{G}$.

Let $\hat{U}\subset\hat{B}$ denote the unipotent radical of $\hat{B}$. The natural morphism $\hat{U}/\hat{B}\to \hat{G}/\hat{G}$ is usually called the Springer resolution.
Let 
\[
S^{\un}_{\hat{G}}=(\hat{U}/\hat{B})\times_{\hat{G}/\hat{G}}(\hat{U}/\hat{B}),
\] 
which we call the (unipotent) Steinberg stack of $\hat{G}$\footnote{As $\hat{U}/\hat{B}\to \hat{G}/\hat{G}$ is not flat, the fiber product needs to considered in derived sense so $S^{\un}_{\hat{G}}$ should be understood as a derived algebraic stack.}. Over $\bC$, there is a $\bG_{m,\bC}$-action on $\hat{U}_\bC$ and therefore on $S^{\un}_{\hat{G},\bC}$, by identifying $\hat{U}_\bC$ with its Lie algebra via the exponential map. Then one can form the quotient stack $S^{\un}_{\hat{G},\bC}/\bG_{m,\bC}$. In the sequel, for an Artin stack $X$ (of finite presentation) over $\bC$, we let $K(X)$ denote the  $K$-group of the ($\infty$-)category of coherent sheaves on $X$.

Kazhdan-Lusztig \cite{Kazhdan-Lusztig-Iwahori} constructed (under the assumption that $G$ is split with connected center) a canonical isomorphism (after choosing a square root of $\sqrt{q}$ of $q$)
\begin{equation}\label{KL isomorphism}
K(S^{\un}_{\hat{G},\bC}/\bG_{m,\bC})\otimes_{K(\bB\bG_{m,\bC})}\bC\cong H^{\cl}_I\otimes\bC
\end{equation}
where the map $K(\bB\bG_{m,\bC})\to \bC$ sends the class corresponding to the tautological representation of $\bG_{m,\bC}$ to $\sqrt{q}$. In addition, the isomorphism induces the Bernstein isomorphism
\begin{equation}\label{Bernstein isomorphism}
K(\bB \hat{G}_\bC)\otimes\bC\cong Z(H^{\cl}_I\otimes\bC)
\end{equation}
where $Z(H^{\cl}_I\otimes\bC)$ is the center of $H^{\cl}_I\otimes\bC$, and the map $K(\bB \hat{G}_\bC)\to K(S^{\un}_{\hat{G},\bC}/\bG_{m,\bC})$ is
induced by the natural projection $S^{\un}_{\hat{G}}/\bG_m\to \bB\hat{G}$.

\begin{rmk}
It would be interesting to give a description of the $\bZ$-algebra $H^{\cl}_I$ in terms of the geometry involving $\hat{G}$, which would generalize the integral Satake isomorphism from \cite{Zhu-Integral-Satake}.
\end{rmk}

It turns out that the Kazhdan-Lusztig isomorphism \eqref{KL isomorphism} also admits a categorification, usually known as the Bezrukavnikov equivalence, which gives two realizations of the affine Hecke category. Again, when switching to the geometric theory, we allow $F$ to be a little bit more general as in Sect. \ref{S: geometric Satake}. We also assume that $G$ extends to a connected reductive group over $\mO$ and fix a pinning of $G$ over $\mO$.
Let $L^+G\to G_k$ be the natural ``reduction mod $\varpi$" map, and let $\Iw\subset L^+G$ be the preimage of $B_k\subset G_k$. 
This is the geometric analogue of $I$. 
Then as in the unramified case discussed in Sect. \ref{S: geometric Satake}, one can define the Iwahori local Hecke stack $\Hk_{\Iw}= \Iw\backslash LG/\Iw$ and the monoidal categories 
$\Shv_{c}(\Hk_{\Iw}\otimes\overline k,\La)\subset\Shv(\Hk_{\Iw}\otimes\overline k,\La)$. The category $\Shv_{c}(\Hk_{\Iw}\otimes\bar k,\La)$ can be thought as a categorical analogue of $H^{\cl}_I$, usually called the affine Hecke category. 

Recall that we let $\breve F=W_\mO(\overline k)[1/\varpi]$. The inertia $I_F:=\Ga_{\breve F}$ of $F$ has a tame quotient $I_F^t$ isomorphic to $ \prod_{\ell\neq p} \bZ_\ell(1)$.
\begin{thm} \label{T: Bezrukavnikov equivalence}
For every choice of a topological generator $\tau$ of the tame inertia $I_F^t$,
there is a canonical equivalence of monoidal $\infty$-categories
\[
\on{Bez}^{\un}_G:   \Coh(S^{\un}_{\hat{G},\bQ_\ell})\cong \Shv_{c}(\Hk_{\Iw}\otimes\overline k,\bQ_\ell).
\]
\end{thm}

In fact, Bezrukavnikov proved such equivalence when $F=k(\!(\varpi)\!)$ in \cite{Bezrukavnikov-Affine-Hecke}. Yun and the author deduce the mixed characteristic case from the equal characteristic case. It would be interesting to know whether the new techniques introduced in \cite{Scholze-Berkeley, Fargues-Scholze-Geometrization} can lead a direct proof of this equivalence in mixed characteristic. (See \cite{AGLR} for some progress in this direction.) 

\begin{rmk}
Again, for arithmetic applications, one needs to describe the action of $\Ga_k$ on $\Shv_{c}(\Hk_{\Iw}\otimes\bar k,\La)$ in terms of the dual group side. See \cite{Bezrukavnikov-Affine-Hecke, Hemo-Zhu-Unipotent} for a discussion. 
\end{rmk}

We explain an important ingredient in the proof of Theorem \ref{T: Bezrukavnikov equivalence} (when $F=k(\!(\varpi)\!)$). There is a smooth affine group scheme $\mG$  (called the Iwahori group scheme) over $\mO$ such that $\mG\otimes F=G$ and $L^+\mG=\Iw$. Then  there is a local Hecke stack $\Hk_{\mG,D}$ over $D$, analogous to $\Hk_{G,D^S}$ as discussed at the end of Sect. \ref{S: geometric Satake} (here $S=\{1\}$). In addition,  $\Hk_{\mG,D}|_{D^*}\cong \Hk_{G,D}|_{D^*}$ and $\Hk_{\mG,D}|_0=\Hk_{\Iw}$, where $0\in D$ is the closed point.
Then taking nearby cycles gives
\begin{equation}\label{E: central sheaf}
\mZ: \Coh(\bB\hat{G}_\La)^\heart\xrightarrow{\Sat_{\{1\}}} \Perv(\Hk_{\mG,D}|_{D^*_{\bar k}},\La)\xrightarrow{\Psi} \Perv(\Hk_{\Iw}\otimes \bar k,\La).
\end{equation} 
This is known as Gaitsgory's central functor \cite{Gaitsgory-Central-Sheaf, Zhu-Coherence-Conjecture}, which can be regarded as a  categorification of \eqref{Bernstein isomorphism}.
We remark this construction is motivated by the Kottwitz conjecture originated from the study of mod $p$ geometry of Shimura varieties. See \S \ref{S: local model} for some discussions. 

Theorem \ref{T: Bezrukavnikov equivalence} admits a generalization to the tame level.  
We consider the following diagram
\[
\hat{G}/\hat{G}\leftarrow \hat{B}/\hat{B}\xrightarrow{q_{\hat{B}}} \hat{T}/\hat{T}.
\]
where the left morphism is the usual Grothendieck-Springer resolution.
Let $\chi$ be a $\La$-point of $\hat{T}/\hat{T}$, where $\La$ is a finite extension of $\bQ_\ell$. Let $(\hat{B}/\hat{B})_\chi=q_{\hat{B}}^{-1}(\chi)$, and let
\[
S^{\chi}_{\hat{G},\La}:= (\hat{B}/\hat{B})_\chi\times_{\hat{G}/\hat{G}} (\hat{B}/\hat{B})_\chi.
\]
Note that if $\chi=1$, this reduces to $S^{\un}_{\hat{G},\La}$.
On the other hand, a (torsion) $\La$-point $\chi\in \hat{T}/\hat{T}$ defines a one-dimensional character sheaf $\mL_\chi$ on $\Iw\otimes\overline{k}$.  Then one can define the monoidal category of bi-$(\Iw,\mL_\chi)$-equivariant constructible sheaves on $LG_{\bar k}$, denoted as $\Shv_{\cons}({}_\chi(\Hk_{\Iw})_{\chi},\La)$. If $\chi=1$ so $\mL_\chi$ is the trivial character sheaf on $\Iw$, this reduces to the affine Hecke category $\Shv_c(\Hk_{\Iw}\otimes\overline k,\La)$.
The following generalization of Theorem \ref{T: Bezrukavnikov equivalence} is conjectured in \cite{Bezrukavnikov-Affine-Hecke} and will be proved in a forthcoming joint work with Dhillon-Li-Yun (\cite{DLYZ-Endoscopy}).

\begin{thm} \label{T: Monodromic Bezrukavnikov equivalence}
Assume that $\cha F=\cha k$. There is a canonical monoidal equivalence
\[
\on{Bez}^{\chi}_{\hat{G}}: \Coh(\hat{S}_{\hat{G},\La}^{\chi})\cong \Shv_{c}({}_\chi(\Hk_{\Iw})_{\chi},\La).
\]
\end{thm}

\begin{rmk}\label{R: int Bez}
It is important to establish a version of equivalences in Theorem \ref{T: Bezrukavnikov equivalence} and \ref{T: Monodromic Bezrukavnikov equivalence} for $\bZ_\ell$-sheaves.
\end{rmk}

\begin{rmk}\label{2-categorical LGL}
The local geometric Langlands correspondence beyond the tame ramification has not been fully understood, although certain wild ramifications have appeared in concrete problems (e.g. \cite{HNY-Kloosterman, Zhu-FG-vs-HNY}).
It is widely believed that the general local geometric Langlands
should be formulated as $2$-categorical statement, predicting the $2$-category of module categories under the action of (appropriately defined) category of sheaves on $LG$ is equivalent to the $2$-category of categories over the stack of local geometric Langlands parameters. The precise formulation is beyond the scope of this surveybut, roughly speaking, it implies (and is more or less equivalent to say) that the Hecke category for appropriately chosen ``level'' of $LG$ is (Morita) equivalent to the category of coherent sheaves on some stack of the form $X\times_{Y}X$, where $Y$ is closely related to the moduli of local geometric Langlands parameters. 
\end{rmk}

\subsection{Global geometric Langlands correspondence}\label{Global GLC}
As mentioned at the beginning of the article, the (global) geometric Langlands program originated from Drinfeld's proof of Langlands conjecture for $\GL_2$ over function fields. Early developments of this subject mostly focused on the construction of Hecke eigensheaves
associated to Galois representations of a global function field $F$ (or equivalently local systems on a smooth algebraic curve $X$), e.g. see \cite{Drinfeld-GL2, Laumon-Geometric-Langlands, FGV-Geometric-Langlands-GLn}. 

The scope of the whole program then shifted after the work \cite{Beilinson-Drinfeld-Hitchin}, in which Beilinson-Drinfeld formulated a rough categorical form of the global geometric Langlands correspondence. The formulation then was made precise by  Arinkin-Gaitsgory in \cite{Arinkin-Gaitsgory-Singular-Support}, which we now recall.
Let $X$ be a smooth projective curve over $F=\bC$.
On the automorphic side, let $\mathbf{D}_c(\Bun_G)$ be the $\infty$-category of coherent D-modules on the moduli stack $\Bun_G$ of principal $G$-bundles on $X$. On the Galois side, let $\Coh(\Loc_{\hat{G}})$ be the $\infty$-category of coherent sheaves on the moduli stack $\Loc_{\hat{G}}$ of de Rham $\hat{G}$-local systems (a.k.a. principal $\hat{G}$-bundles with flat connection) on $X$.

\begin{conjecture}\label{Conj: GGLC}
There is a canonical equivalence of $\infty$-categories
\[
\bL_G:  \Coh(\Loc_{\hat{G}})\cong \mathbf{D}_c(\Bun_G),
\]
satisfying a list of natural compatibility conditions.
\end{conjecture}
We briefly mention the most important compatibility condition, and refer to \cite{Arinkin-Gaitsgory-Singular-Support} for the rest. Note that both sides admit actions by a family of commuting operators labelled by $x\in X$ and $V\in \Coh(\bB\hat{G}_\bC)^\heart$. Namely,
for every point $x\in X$, there is the evaluation map $\Loc_{\hat{G}}\to \bB\hat{G}_\bC$ so every 
$V\in \Coh(\bB\hat{G}_\bC)^\heart$ gives a vector bundle $\widetilde{V}_x$ on $\Loc_{\hat{G}}$ by pullback, which then acts on $\Coh(\Loc_{\hat{G}})$ by tensoring. On the other hand, there is the Hecke operator $H_{V,x}$ that acts on $\mathbf{D}_c(\Bun_G)$ by convolving the sheaf $\Sat_{\{1\}}(V)|_x$ from the ($D$-module version of) the geometric Satake \eqref{factorizable Sat}. Then the equivalence $\bL_G$ should intertwine the actions of these operators.

Although the conjecture remains widely open, it is known that the category of perfect complexes $\Perf(\Loc_{\hat{G}})$ on $\Loc_{\hat{G}}$ acts on $\mathbf{D}_c(\Bun_G)$, usually called the spectral action, in a way such that the action of $\widetilde{V}_x\in \Perf(\Loc_{\hat{G}})$ on $\mathbf{D}_c(\Bun_G)$ is given by the Hecke operator $H_{V,x}$.

Nowadays, Conjecture \ref{Conj: GGLC} sometimes is referred as the de Rham version of the global geometric Langlands conjecture, as there are some other versions of such conjectural equivalences, which we briefly mention.

First, in spirit of the non-abelian Hodge theory, there should exist a classical limit of Conjecture \ref{Conj: GGLC}, sometimes known as the Dolbeault version of the global geometric Langlands. While the precise formulation is unknown (to the author), generically, it amounts to the duality of Hitchin fibrations for $G$ and $\hat{G}$ (in the sense of mirror symmetry), and was established ``generically" in \cite{Donagi-Pantev-Hitchin, Chen-Zhu-Geometric-Langlands-char-p}. By twisting/deforming such duality in positive characteristic, one can prove a characteristic $p$ analogue of Conjecture \ref{Conj: GGLC} (of course only  ``generically", see
\cite{Bezrukavnikov-Braverman-Geometric-Langlands-char-p, Chen-Zhu-Nonabelian-Hodge, Chen-Zhu-Geometric-Langlands-char-p}). Interestingly, this ``generic" characteristic $p$ equivalence can be used to give a new proof of the main result of  \cite{Beilinson-Drinfeld-Hitchin} (at least for $G=\GL_n$, see \cite{Bezrukavnikov-Travkin-Quantization}).

The work \cite{Beilinson-Drinfeld-Hitchin} (and therefore the de Rham version of the global geometric Langlands) was strongly influenced by conformal field theory. On the other hand,
motivated by topological field theory, Ben-Zvi and Nadler \cite{Ben-Zvi-Nadler-Betti-Geometric-Langlands} proposed a Betti version of Conjecture \ref{Conj: GGLC}, where on the automorphic side the category of $D$-modules on $\Bun_G$ is replaced with the category of sheaves of $\bC$-vector spaces on (the analytification of) $\Bun_G$ and on the Galois side $\Loc_{\hat{G}}$ is replaced by the moduli of Betti $\hat{G}$-local systems (a.k.a. $\hat{G}$-valued representations of fundamental group of $X$). 

The Riemann-Hilbert correspondence allows one to pass between the de Rham $\hat{G}$-local systems and Betti $\hat{G}$-local systems, but in a transcendental way. So Conjecture \ref{Conj: GGLC} and its Betti analogue are not directly related.
Recently, Arinkin et. al. \cite{AGKRRV} proposed yet another version of Conjecture \ref{Conj: GGLC}, which directly relates both de Rham and Betti version, and at the same time includes a version in terms of $\ell$-adic sheaves. So it is more closely related to the classical Langlands correspondence over function fields, as will be discussed in Sect. \ref{S: Global Langlands function fields}.

\section{From geometric to classical Langlands program}
In the previous section, we discussed how the ideas of categorification and geometrization lead to developments of the geometric Langlands program. On the other hand,
the ideas of quantum physics allow one to reverse arrows in \eqref{categorification}
by evaluating a (topological) quantum field theory at manifolds of different dimensions. Such ideas are certainly not new in geometry and topology. But surprisingly, it also leads to new understanding of the classical Langlands program. Indeed, it has been widely known that there is analogy between global fields and $3$-manifolds, and under such analogy Frobenius corresponds to the fundamental group of a circle. Then ``compactification of field theories on a circle" leads to the categorical trace method (e.g. see \cite{AGKRRV, BCHN, Zhu-CDM}), which plays more and more important roles in geometric representation theory.

\subsection{Categorical arithmetic local Langlands}\label{S: Local Langlands}

In this subsection, let $F$ be either a finite extension of $\bQ_p$ or is isomorphic to $\bF_q((\varpi))$.  The classical local Langlands correspondence seeks a classification of smooth irreducible representations of $G(F)$ in terms of Galois data. The precise formulation, beyond the $G=\GL_n$ case (which is a theorem by \cite{LRS-GLn, Harris-Taylor-GLn}), is complicated.
However, the yoga that the local geometric Langlands is $2$-categorical (see Remark \ref{2-categorical LGL}) suggests that the even the classical local Langlands correspondence should and probably needs to be categorified.

The first ingredient needed to formulate the categorical arithmetic local Langlands is the following result, due independently by \cite{DHKM, Fargues-Scholze-Geometrization, Zhu-Conjecture}. We take the formulation from \cite{Zhu-Conjecture}, and refer to \emph{loc. cit.} for the notion of (strongly) continuous homomorphisms.

\begin{thm}\label{Stack of local parameters}
The prestack sending a $\bZ_\ell$-algebra $A$ to the space of (strongly) continuous homomorphisms $\rho: W_F\to \CG(A)$ such that $d\circ\rho=(\cycl^{-1},\pr)$ is represented by a (classical) scheme $\Loc_{\CG}^\Box$, which is a disjoint union of affine schemes that are flat, of finite type, and of locally complete intersection over $\bZ_\ell$. 
\end{thm}

The conjugation action of $\hat{G}$ on $\CG$ induces an action of $\hat{G}$ on $\Loc_{\CG}^\Box$ and we call the quotient stack $\Loc_{\CG}=\Loc_{\CG}^\Box/\hat{G}$ be the stack of local Langlands parameters, which roughly speaking classifies the groupoid of the above $\rho$'s up to $\hat{G}$-conjugacy. 

In the categorical version of the local Langlands correspondence,
on the Galois side it is natural to consider the ($\infty$-)category $\Coh(\Loc_{\CG})$ of coherent sheaves on $\Loc_{\CG}$.
On the representation side, one might naively consider the ($\infty$-)category $\Rep(G(F),\La)$ of smooth representations of $G(F)$.
But in fact, this category needs to be enlarged.
This can be seen from different point of views. Indeed, it is a general wisdom shared by many people that in the classical local Langlands correspondence, it is better to study representations of $G$ together with a collection of groups that are (refined version of) its inner forms. There are various proposals of such collections.
Arithmetic geometry (i.e. the study of $p$-adic and mod $p$ geometry of Shimura varieties and moduli of Shtukas) and geometric representation theory (i.e. the categorical trace construction) suggest  to study a category glued from the categories of representations of  a collection of groups $\{J_b(F)\}_{b\in B(G)}$ arising from the Kottwitz set 
\[
B(G) = G(\breve F)/\sim, \quad \mbox{where } g\sim g' \mbox{ if } g'=h^{-1}g\sigma(h) \mbox{ for some }  h\in G(\breve F).
\] 
Here for $b\in B(G)$ (lifted to an element in $G(\breve F)$), the group $J_b$  is an $F$-group defined by assigning and $F$-algebra the group $J_b(R)=\{h\in G(\breve F\otimes_FR)\mid h^{-1}b\sigma(h)=b\}$. In particular, if $b=1$ then $J_b=G$. In general, there is a well-defined embedding $(J_{b})_{\overline{F}}\to G_{\overline{F}}$ up to conjugacy, making $J_b$ a refinement of an inner form of a Levi subgroup of $G$ (say $G$ is quasi-split).

There are two ways to make this idea precise.
One is due to Fargues-Scholze \cite{Fargues-Scholze-Geometrization}, who regard $B(G)$ as the set of points of the $v$-stack $\Bun_G$ of $G$-bundles on the Fargues-Fontaine curve and consider the category $D_{\on{lis}}(\Bun_G,\La)$ of appropriately defined \'etale sheaves on $\Bun_G$, which indeed glues all $\Rep(J_b(F),\La)$'s together. We mention that this approach relies on Scholze's work on $\ell$-adic formalism of diamond and condensed mathematics. 

In another approach \cite{Xiao-Zhu-Cycle, Zhu-CDM, Zhu-Conjecture, Hemo-Zhu-Unipotent}, closely related to the idea of categorical trace, the set $B(G)$ is regarded as the set of points of the (\'etale) quotient stack 
\begin{equation*}\label{E: kot stack}
\frakB(G):=LG/\Ad_\sigma LG,
\end{equation*} 
where $\Ad_\sigma$ denotes the Frobenius twisted conjugation given by 
\[
\Ad_\sigma: LG\times LG\to LG,\quad  (h,g)\mapsto hg\sigma(h)^{-1}.
\]
Then we have the category of $\La$-sheaves $\Shv(\kot(G)\otimes \bar{k},\La)$ as mentioned before. Although $\kot(G)$ is a wild object in the traditional algebraic geometry, there are still a few things one can say about its geometry and
the category $\Shv(\frakB(G)\otimes \bar{k},\La)$ is quite reasonable. In addition, it is possible to define the category $\Shv_c(\frakB(G)\otimes \bar{k},\La)$ of constructible sheaves on $\frakB(G)\otimes \bar{k}$, as we now briefly explain and refer to \cite{Hemo-Zhu-Unipotent} for careful discussions.

For every algebraically closed field $\Omega$ over $k$, the groupoid of $\Omega$-points of $\kot(G)$ is the groupoid of $F$-isocrystals with $G$-structure over $\Omega$ and the set of its isomorphism classes can be identified with the Kottwitz set $B(G)$. However, $\kot(G)$ is not merely a disjoint union of its points. Rather, it admits a stratification, known as the Newton stratification, labelled by $B(G)$. Namely, the set $B(G)$ has a natural partial order and roughly speaking for each $b\in B(G)$ those $\Omega$-points corresponding to $b'\leq b$ form a closed substack $i_{\leq b}: \kot(G)_{\leq b}\subset \kot(G)\otimes \bar k$ and those $\Omega$-points corresponding to $b$ form an open substack $j_{b}: \kot(G)_b\subset \kot(G)_{\leq b}$.  
In particular, basic elements in $B(G)$ (i.e. minimal elements with respect to the partial order $\leq$) give closed strata. We also mention that if $b$ is basic, $J_b$ is a refinement of an inner form of $G$, usually called an extended pure inner form of $G$.

In the rest of this subsection, we simply write $\kot(G)\otimes \bar k$ by $\kot(G)$. We write $i_b=i_{\leq b}j_b: \kot(G)_b\hookrightarrow \kot(G)$ for the locally closed embedding. For $b$, let $\fgRep(J_b(F),\La)$ be the full subcategory of $\Rep(J_b(F),\La)$ generated (under finite colimits and retracts) by compactly induced representations 
\[
\delta_{K,\La}:=c\on{-ind}_K^{J_b(F)}(\La)
\]
from the trivial representation of open compact subgroups $K\subset J_b(F)$.   
The following theorem from \cite{Hemo-Zhu-Unipotent} summarizes some properties of $\Shv_c(\kot(G),\La)$. 

\begin{thm}\label{T: Local-Langlands-Category}
\begin{enumerate}
    \item An object  in $\Shv(\kot(G),\La)$ is constructible if and only if its $!$-restriction to each $\frakB(G)_b$ is constructible and is zero for almost all $b$'s. If $\La$ is a field of characteristic zero, $\Shv_{c}(\frakB(G),\La)$ consist of compact objects in $\Shv(\kot(G),\La)$. 
    \item For every $b\in B(G)$, there is a canonical equivalence 
    $\Shv_c(\frakB(G)_b,\La)\cong \fgRep(J_b(F),\La)$. There are fully faithful embeddings $i_{b,*}, i_{b,!}:  \Shv_c(\frakB(G)_b,\La)\to \Shv_c(\kot(G),\La)$  (which coincide when $b$ is basic), inducing a semi-orthogonal decomposition of $\Shv_c(\frakB(G),\La)$ in terms of $\{\Shv_c(\frakB(G)_b,\La)\}_b$.     
    \item There is a self-duality functor
    $\bD^{\coh}\colon \Shv_c(\frakB(G),\La)\simeq \Shv_c(\frakB(G),\La)^{\vee}$ obtained by gluing cohomological dualities (in the sense of Bernstein-Zelevinsky) on various $\fgRep(J_b(F),\La)$'s.
    \item There is a natural perverse $t$-structure obtained by gluing (shifted) $t$-structures on various $\fgRep(J_b(F),\La)$'s, preserved by $\bD^{\coh}$ if $\La$ is a field. 
\end{enumerate} 
\end{thm}

The following categorical form of the arithmetic local Langlands conjecture (\cite[Sect. 4.6]{Zhu-Conjecture}) is inspired by the global geometric Langlands conjecture as discussed in \S \ref{Global GLC}.

\begin{conjecture}\label{C: LLC}
Assume that $G$ is quasi-split over $F$ equipped with a pinning $(B,T,e)$ and fix a non-trivial additive character $\psi: F\to \bZ_\ell[\mu_{p^\infty}]^\times$.
There is a canonical equivalence of categories
\[
\bL_{G}:   \Coh(\Loc_{\CG}\otimes\La)\cong \Shv_{c}(\frakB(G),\La). 
\]
\end{conjecture}

\begin{rmk}
(1) There is a closely related version of the conjecture, with $\Shv_{c}(\frakB(G),\La)$ replaced by $\Shv(\frakB(G),\La)$ and with $\Coh(\Loc_{\CG}\otimes\La)$ replaced by its ind-completion (with certain support condition imposed) (see \cite[Sect. 4.6]{Zhu-Conjecture}).
Fargues-Scholze \cite{Fargues-Scholze-Geometrization} make a conjecture parallel to this version, with the category $\Shv(\frakB(G),\La)$ replaced by $D_{\on{lis}}(\Bun_G,\La)$ as mentioned above.

(2) It is also explained in \cite{Zhu-Conjecture} a motivic hope to have a version of such equivalence over $\bQ$.
\end{rmk}
One consequence of the conjecture is that for every $b$, there should exist a fully faithful embedding
\begin{equation*}\label{E: functor A}
\frakA_{J_b}:  \fgRep(J_b(F),\La)\to  \Coh(\Loc_{\CG}\otimes\La),
\end{equation*}
obtained as the restriction of a quasi-inverse of $\bL_G$ to $i_{b,!}(\fgRep(J_b(F),\La))$.
The existence of such functor is closely related to the idea of local Langlands in families and has also been considered (in the case $J_b=G$ is split and $\La$ is a field of characteristic zero) in \cite{Hellmann-Iwahori, BCHN}. 

In particular, for every open compact subgroup $K\subset J_b(F)$ there should exist a coherent sheaf 
\begin{equation}\label{E: sheaf A}
\frakA_{K,\La}:=\frakA_{J_b}(\delta_{K,\La})
\end{equation} 
on $\Loc_{\CG}\otimes\La$, such that
\begin{equation}\label{E: Hecke alg isom}
\big(R\End_{\Coh(\Loc_{\CG}\otimes\La)}\frakA_{K,\La}\big)^{\op}\cong \big(R\End_{\Rep(G(F),\La)}\delta_{K,\La}\big)^{\op}=:H_{K,\La}.
\end{equation}
The algebra $H_{K,\La}$ is sometimes called the derived Hecke algebra as it might not concentrate on cohomological degree zero (when $\La=\bZ_\ell$ or $\bF_\ell$). 
See \cite[Sect. 4.3-4.5]{Zhu-Conjecture} for conjectural descriptions of $\frakA_{K,\La}$ in various cases.

As in the global geometric Langlands conjecture, the equivalence from Conjecture \ref{C: LLC} should satisfy a set of compatibility conditions. 
For example, it should be compatible with parabolic inductions on both sides, and should be compatible with the duality $\bD^{\coh}$ on $\Shv_c(\frakB(G),\La)$ and the (modified) Grothendieck-Serre duality of $\Coh(\Loc_{\CG}\otimes\La)$.  We refer to \cite{Zhu-Conjecture, Hemo-Zhu-Unipotent} for more details.

On the other hand, Conjecture \ref{C: LLC} predicts an action of the category $\Perf(\Loc_{\CG}\otimes\La)$ of perfect complexes on $\Loc_{\CG}\otimes\La$ on $\Shv_c(\frakB(G),\La)$, analogous to the spectral action as mentioned in Sect. \ref{Global GLC}. One of the main results of \cite{Fargues-Scholze-Geometrization} is the construction of such action in their setting. 
Currently the existence of such spectral action on $\Shv_c(\frakB(G),\La)$ is not known. But there are convincing evidences that
Conjecture \ref{C: LLC} should still be true.

We assume that $G$ extends to a reductive group over $\mO$ as before. Then there are closed substacks 
\[\Loc_{\CG}^{\ur}\subset\Loc_{\CG}^{\un}\subset \Loc_{\CG},
\] 
usually called the stack of unramified parameters (resp. unipotent parameters), classifying those $\rho$ such that $\rho(I_F)$ is trivial (resp. $\rho(I_F)$ is unipotent).
For $\La=\bQ_\ell$, $\Loc_{\CG}^{\un}\otimes\bQ_\ell$ is a connected component of $\on{Loc}_{\CG}\otimes\bQ_\ell$. 

On the other hand, there is
the unipotent subcategory $\Shv_c^{\un}(\kot(G),\bQ_\ell)\subset \Shv_c(\kot(G),\bQ_\ell)$, which roughly speaking is the glue of categories $\fgRep^{\un}(J_b(F),\bQ_\ell)$ of unipotent representations of $J_b(F)$ (introduced in \cite{Lusztig-Unipotent}) for all $b\in B(G)$. We have the following theorem from \cite{Hemo-Zhu-Unipotent}, deduced from Theorem \ref{T: Bezrukavnikov equivalence} by taking the Frobenius-twisted categorical trace. 

\begin{thm}\label{T: tame local Langlands}
For a reductive group $G$ over $\mO$ with a fixed pinning $(B,T,e)$, there is a canonical equivalence
\[
\bL^{\un}_G:  \Coh(\Loc_{\CG}^{\un}\otimes\bQ_\ell)\cong  \Shv_{c}^{\un}(\frakB(G),\bQ_\ell).
\]
\end{thm}
For arithmetic applications, it is important to match specific objects under the equivalence. We give a few examples and refer to  \cite{Hemo-Zhu-Unipotent} for many more of such matchings (see also \cite[Sect. 4.3-4.5]{Zhu-Conjecture}).

\begin{ex}
The equivalence $\bL_G^{\un}$ gives the the conjectural coherent sheaf in \eqref{E: sheaf A} for all parahoric subgroups $K\subset G(F)$ (in the sense of Bruhat-Tits) such that \eqref{E: Hecke alg isom} holds. For example, we have $\frakA_{G(\mO),\bQ_\ell}\cong \mO_{\Loc_{\CG}^{\ur}}\otimes\bQ_\ell$, which gives
\begin{equation}\label{E: derived Satake}
\big(R\End_{\Coh(\Loc_{\CG})}\mO_{\Loc_{\CG}^{\ur}}\big)^{\op}\otimes\bQ_\ell\cong \big(R\End\ \delta_{G(\mO),\bQ_\ell}\big)^{\op}= H_{G(\mO),\bQ_\ell}\cong H_{G(\mO)}^{\cl}\otimes\bQ_\ell.
\end{equation} 
As $\Loc_{\CG}^{\ur}\cong (\CG|_{d=(q,\sigma)})/\hat{G}$,
taking $0$th cohomology recovers the Satake isomorphism \eqref{Sat isom}. 
In addition, it implies that the left hand side has no higher cohomology, which is not obvious.  We mention that it is conjectured in \cite[Sect. 4.3]{Zhu-Conjecture} that $\frakA_{G(\mO),\bZ_\ell}\cong \mO_{\Loc_{\CG}^{\ur}}$ so the first isomorphism in \eqref{E: derived Satake} should hold over $\bZ_\ell$, known as the (conjectural) derived Satake isomorphism. (But $H_{G(\mO),\bZ_\ell}\neq H_{G(\mO)}^{\cl}\otimes\bZ_\ell$ in general.)

There is also a pure Galois side description of $\frakA_{I,\bQ_\ell}$, known as the unipotent coherent Springer sheaf as defined in \cite{BCHN, Zhu-Conjecture} (see also \cite{Hellmann-Iwahori}).
\end{ex}

\begin{ex}\label{Ex: trace of Satake}
By construction, there is a natural morphism of stacks $\Loc_{\CG}\to \bB\hat{G}$ over $\bZ_\ell$. For a representation of $\hat{G}$ on a finite projective $\La$-module, regarded as a vector bundle on $\bB\hat{G}_\La$, let $\widetilde V$ be its pullback to $\Loc_{\CG}\otimes\La$, and let $\widetilde{V}^{?}\in \Perf(\Loc_{\CG}^{?}\otimes\La)$ be its restriction of $\Loc_{\CG}^?\otimes\La$ for $?=\ur$ or $\un$. 
Note that for $\La=\bQ_\ell$, $\widetilde{V}^{\ur}\cong \widetilde{V}\otimes \frakA_{G(\mO),\bQ_\ell}$. We have
\[
\bL_G^{\un}(\widetilde{V}^{\ur})\cong \on{Nt}_!r^!\Sat(V)=:\mS_V
\]
where $r$ and $\on{Nt}$ are maps in the following correspondence.
\[
\Hk_G=L^+G\backslash LG/L^+G\xleftarrow{r} LG/\Ad_\sigma L^+G\xrightarrow{\on{Nt}} LG/\Ad_\sigma LG=\kot(G)
\]
In particular, for two representations $V$ and $W$ of $\hat{G}$, there is a morphism
\begin{equation}\label{E: local S-operator}
R\Hom_{\Loc_{\CG}^{\ur}\otimes\bQ_\ell}(\widetilde{V}^{\ur}, \widetilde{W}^{\ur})\to R\Hom_{\Shv_c(\kot(G),\bQ_\ell})(\mS_V,\mS_W)
\end{equation}
compatible with compositions. Such map was first constructed in \cite{Xiao-Zhu-Cycle,Zhu-CDM} and (the version for underived $\Hom$ spaces) was then extended to $\bZ_\ell$-coefficient in \cite{Yu-Jacquet-Langlands-Higher-Weights}. It has significant arithmetic applications, as will be explained in Sect. \ref{Application}. 
\end{ex}

\begin{rmk}
Likely Theorem \ref{T: tame local Langlands} can be extended to the tame level by  taking the Frobenius-twisted categorical trace of the equivalence from Theorem \ref{T: Monodromic Bezrukavnikov equivalence}. On the other hand, as mentioned in Remark \ref{R: int Bez}, it is important to extend these equivalences to $\bZ_\ell$-coefficient.
\end{rmk}

\subsection{Global arithmetic Langlands for function fields}\label{S: Global Langlands function fields}
Next we turn to global aspects of the arithmetic Langlands correspondence. As mentioned at the beginning, its classical formulation very roughly speaking predicts a natural correspondence between the set of (irreducible) Galois representations and the set of (cuspidal) automorphic representations. As in the local case, beyond the $\GL_n$ case (which is a theorem by \cite{Lafforgue-GLn}), such formulation is not easy to be made precise. On the other hand,
the global geometric Langlands conjecture from Sect. \ref{Global GLC} and philosophy of decategorification/trace suggest that the global arithmetic Langlands can and probably should be formulated as an isomorphism between two vector spaces, arising from the Galois and the automorphic side respectively. In this subsection, we formulate such a conjecture in the global function field case.

Let $F=\bF_q(X)$ be the function field of a geometrically connected smooth projective curve $X$ over $\bF_q$.
We write $\eta=\Spec F$ for the generic point of $X$ and $\overline\eta$ for a geometric point over $\eta$. Let $|X|$ denote the set of closed points of $X$. For $v\in |X|$, let $\mO_v$ denote the complete local ring of $X$ at $v$ and $F_v$ its fractional field. Let
$\bO_F=\prod_{v\in|X|}\mO_v$ be the integral ad\`eles, and $\bA_F=\prod'_{v\in|X|}F_v$ the ring of ad\`eles. For a finite non-empty set of places $Q$, let $W_{F,Q}$ denote the Weil group of $F$, unramified outside $Q$. 

Let $G$ be a connected reductive group over $F$. Similarly to the local situation, the first step to formulate our global conjecture is the following theorem from \cite{Zhu-Conjecture}.
\begin{thm}
Assume that $\ell\nmid 2p$.  The prestack sending a $\bZ_\ell$-algebra $A$ to the space of (strongly) continuous homomorphisms $\rho: W_{F,Q}\to \CG(A)$ such that $d\circ\rho=(\cycl^{-1},\pr)$ is represented by a derived scheme $\Loc_{\CG,Q}^\Box$, which is a disjoint union of derived affine schemes that are flat and of finite type over $\bZ_\ell$. If $Q\neq\emptyset$, $\Loc_{\CG,Q}^\Box$ is quasi-smooth. 
\end{thm}

We then define the stack of global Langlands parameters as $\Loc_{\CG,Q}=\Loc_{\CG,Q}^\Box/\hat{G}$. Similar to the local case (see Example \ref{Ex: trace of Satake}),
for a representation of $\hat{G}_\La$ on a finite projective $\La$-module, regarded as a vector bundle on $\bB\hat{G}_\La$, let $\widetilde V$ be its pullback to $\Loc_{\CG,Q}\otimes\La$. If $V$ is the restriction of a representation of $(\CG)^S$ along the diagonal embedding $\hat{G}\to (\CG)^S$, then there is a natural (strongly) continuous $W_{F,Q}^S$-action on $\widetilde V$ (see \cite[Sect. 2.4]{Zhu-Conjecture}). For a place $v$ of $F$, let $\Loc_{\CG,v}$ denote the stack of local Langlands parameters for $G_{F_v}$. Let
\[
\res: \Loc_{\CG,Q}\to \prod_{v\in Q}\Loc_{\CG,v}
\] 
denote the map by restricting global parameters to local parameters (induced by the map $W_{F_v}\to W_{F,Q}$). Later on, we will consider the $!$-pullback of coherent sheaves on $ \prod_{v\in Q}\Loc_{\CG,v}$ along this map.
\begin{rmk}
(1) In fact, when $Q=\emptyset$, the definition of $\Loc_{\CG,Q}$ needs to be slightly modified. 

(2) Unlike the local situation, $\Loc_{\CG,Q}$ has non-trivial derived structure in general (see \cite[Rem. 3.4.5]{Zhu-Conjecture}). Let ${}^{cl}\Loc_{\CG,Q}$ denote the underlying classical stack.

(3) A different definition of $\Loc_{\CG,Q}\otimes\bQ_\ell$ is given by \cite{AGKRRV}.
\end{rmk}

Next we move to the automorphic side.
For simplicity, we assume that $G$ is split over $\bF_q$ in this subsection. 
Fix a level, i.e. an open compact subgroup $K\subset G(\bO_F)$. Let $Q$ be the set of places consisting of those $v$ such that $K_v\neq G(\mO_v)$.
For a finite set $S$, let $\Sht_K(G)_{(X-Q)^S}$ denote the ind-Deligne-Mumford stack over $(X-Q)^S$ of the moduli of $G$-shtukas on $X$ with $S$-legs in $X-Q$ and $K$-level structure. (E.g. see \cite{Lafforgue-Automorphic-to-Galois} for basic constructions and properties of this moduli space.)
Its base change along the diagonal map $\overline\eta\to (X-Q)\xrightarrow{\Delta}(X-Q)^S$ is denoted by $\Sht_K(G)_{\Delta(\overline\eta)}$. 
For every representation $V$ of $(\CG)^S$ on a finite projective $\La$-module, the geometric Satake \eqref{factorizable Sat} (with $D$ replaced by $X-Q$ and with $\La=\bZ_\ell$ allowed) provides a perverse sheaf $\on{Sat}_S(V)$ on $\Sht_K(G)_{(X-Q)^S}$. 
Let $C_c\bigl(\Sht_K(G)_{\Delta(\overline\eta)},\Sat_S(V)\bigr)$ denote the (cochain complex of the) total compactly supported cohomology of $\Sht_K(G)_{\Delta(\overline\eta)}$ with coefficient in $\Sat_S(V)$. It admits a (strongly) continuous action of $W_{F,Q}^S$ (see \cite{HRS-Kunneth} for the construction of such action at the derived level, based on \cite{Xue-Finiteness, Xue-Smoothness}), as well as an action of the corresponding global (derived) Hecke algebra (with coefficients in $\La$) 
\begin{equation}\label{E: derived hecke global}
H_{K,\La}=\Big(R\End \big(c\on{-ind}_K^{G(\bA_F)}(\La)\big)\Big)^{\op}.
\end{equation}
For example,
if $V=\mathbf{1}$ is the trivial representation, then (under our assumption that $G$ is split)
$$C_c\bigl(\Sht_K(G)_{\Delta(\overline\eta)}, \on{Sat}_{\{1\}}(\mathbf{1})\bigr)= C_c(G(F)\backslash G(\bA)/K,\La).$$
Here $G(F)\backslash G(\bA)/K$ is regarded as a discrete DM stack over $\overline\eta$, and $C_c(G(F)\backslash G(\bA)/K,\La)$ denotes its compactly supported cohomology.
When $\La=\bQ_\ell$, this is the space of compactly supported functions on $G(F)\backslash G(\bA)/K$. 

We will fix a pinning $(B,T,e)$ of $G$ and a non-degenerate character $\psi: F\backslash \bA\to \bZ_\ell[\mu_{p}]^\times$, which gives the conjectural equivalence $\bL_{v}$ as in Conjecture \ref{C: LLC} for every $v\in Q$.  In particular, corresponding to $K_v\subset G(F_v)$ there is a conjectural coherent sheaf $\frakA_{K_v}$ (see \eqref{E: sheaf A}) on $\Loc_{\CG,v}$. 

\begin{conjecture}\label{C: gen coh}
There is a natural $(W_{F,Q}^S\times H_{K,\La})$-equivariant isomorphism
\[
R\Gamma\bigl(\on{Loc}_{\CG,Q}\otimes\La, \widetilde{V}\otimes  \res^!(\boxtimes_{v\in Q} \frakA_{K_v})\bigr)\cong C_c\bigl(\on{Sht}_K(G)_{\Delta(\overline\eta)}, \on{Sat}_S(V)\bigr).
\]
\end{conjecture}

We refer to \cite[Sect. 4.7]{Zhu-Conjecture} for more general form of the conjecture (where ``generalized level structures" are allowed) and examples of such conjecture in various special cases. This conjecture could be regarded a precise form of the global Langlands correspondence for function fields. Namely, it gives a precise recipe to match Galois representations and automorphic representations. (E.g. V. Lafforgue's excursion operators are encoded in such isomorphism, see below.)
Moreover, such isomorphism fits in the Arthur-Kottwitz multiplicity formula and at the same time extends such formula to the integral level and therefore relates to automorphic lifting theories. 

The most appealing evidence of this conjecture is the following theorem \cite{Lafforgue-Zhu-Elliptic, Zhu-Conjecture}, as suggested (at the heuristic level) by Drinfeld as an interpretation of Lafforgue's construction. 

\begin{thm}\label{T: gen LZ}
For each $i$, there is a quasi-coherent sheaf $\frakA^i_K$ on ${}^{cl}\Loc_{\CG,Q}\otimes\bQ_\ell$, equipped with an action of $H_{K,\bQ_\ell}$, such that for every finite dimensional $\bQ_\ell$-representation $V$ of $(\CG)^S$, there is a natural $(W_{F,Q}^S\times H_{K,\bQ_\ell})$-equivariant isomorphism
\begin{equation*}\label{E: coh Sht}
\Gamma\big({}^{cl}\Loc_{{}^cG,Q}\otimes\bQ_\ell,\widetilde V\otimes \frakA^i_K\big)\cong H^i_c\bigl(\Sht_K(G)_{\Delta(\overline\eta)},\Sat_S(V)\bigr).
\end{equation*}
\end{thm}
We mention that this theorem actually was proved for any $G$ in \cite{Lafforgue-Zhu-Elliptic, Zhu-Conjecture}. In addition, when $K$ is everywhere hyperspecial, \eqref{E: coh Sht} holds at the derived level by \cite{AGKRRV}.
 
The isomorphism \eqref{E: coh Sht} induces an action of $\Gamma\big({}^{cl}\Loc_{{}^cG,Q}\otimes\bQ_\ell,\mO\big)$ on the right hand side. This is exactly the action by V. Lafforgue's excursion operators, which induces the decomposition of the right hand side (in particular $C_c(G(F)\backslash G(\bA_F)/K, \bQ_\ell)$) in terms of semisimple Langlands parameters. As explained \cite{Lafforgue-Zhu-Elliptic}, over an elliptic Langlands parameter, such isomorphism is closely related to the Arthur-Kottwitz multiplicity formula.  In the case of $G=\GL_n$,  it gives the following corollary, generalizing \cite{Lafforgue-GLn}.
\begin{cor}
Let $\pi$ be a cuspidal automorphic representation of $\GL_n$, with the associated irreducible Galois representation $\rho_\pi: W_{F,Q}\to \GL_n(\La)$ for some finite extension $\La/\bQ_\ell$ and with $\frakm_\pi$ the corresponding maximal ideal of $\Gamma\big({}^{cl}\Loc_{{}^cG,Q}\otimes\La,\mO\big)$. Then there is an $(W_{F,Q}^S\times H_K)$-equivariant isomorphism
\[
H^*_c\big(\Sht_K(G)_{\Delta(\overline\eta)}, \Sat_S(V)\big)/\frakm_\pi\cong  V_{\rho_\pi}\otimes \pi^K.
\]
In particular, the left hand side only concentrates in cohomological degree zero.
\end{cor}

\subsection{Geometric realization of Jacquet-Langlands transfer}\label{geometric JL}
The global Langlands correspondence for number fields is far more complicated. In fact, there are analytic part of the theory which currently seems not to fit the categorification/decategorification framework. Even if we just restrict to the algebraic/arithmetic part of the theory, there are complications coming from the place at $\ell$ and at $\infty$. In particular, the categorical forms of the local Langlands correspondence at $\ell$ and $\infty$ are not yet fully understood.

Nevertheless, in a forthcoming joint work with Emerton and Emerton-Gee (\cite{Emerton-Zhu, Emerton-Gee-Zhu}), we will formulate conjectural Galois theoretical descriptions for the cohomology of Shimura varieties and even cohomology for general locally symmetric space, parallel to Conjecture \ref{C: gen coh}. In this subsection, we just review a conjecture from \cite{Zhu-Conjecture} on the geometric realization of Jacquet-Langlands transfer via cohomology of Shimura varieties and discuss results from \cite{Xiao-Zhu-Cycle, Hemo-Zhu-Unipotent} towards this conjecture. 

We fix a few notations and assumptions. We fix a prime $p$ in this subsection.
Let $\bA_f=\prod'_q\bQ_q$ denote the ring of finite ad\`eles of $\bQ$, and $\bA_f^p=\prod'_{q\neq p}\bQ_p$.  We write $\overline\eta=\Spec \overline{\bQ}$, where $\overline{\bQ}$ is the algebraic closure of $\bQ$ in $\bC$. For a Shimura datum $(G,X)$, let $\mu$ be the (minuscule) dominant weight of $\hat{G}$ (with respect to $(\hat{B},\hat{T})$) determined by $(G,X)$ in the usual way and let $V_\mu$ denote the minuscule representation of $\hat{G}$ of highest weight $\mu$. Let $E\subset\overline{\bQ}\subset\bC$ be the reflex field of $(G,X)$ and write $d_\mu=\dim X$.
For a level (i.e. an open compact subgroup) $K=K_pK^p\subset G(\bQ_p)G(\bA_f^p)$,
let $\Sh_K(G)$ be the corresponding Shimura variety of level $K$ (defined over the reflex field $E$), and let $\Sh_K(G)_{\overline\eta}$ denote its base change along $E\to \overline\bQ$.
Let $v$ be a place of $E$ above $p$. By a specialization $\on{sp}: \overline\eta\to \overline{v}$, we mean a morphism from $\overline\eta$ to the strict henselianization of $\mO_E$ at $v$.
 
To avoid many complications from Galois cohomology (e.g. the difference between extended pure inner forms and inner forms) 
and also some complications from geometry (e.g. the relation between Shimura varieties and moduli of Shtukas), we assume that $G$ is of adjoint type in the rest of this subsection, and refer to \cite{Xiao-Zhu-Cycle} for general $G$. See also \cite{Zhu-Conjecture} with less restrictions on $G$.

\begin{dfn}\label{D: prime-to-p trivial}
Let $G$ be a connected reductive group over $\bQ$. A prime-to-$p$ (resp. finitely) trivialized inner form of $G$ is a $G$-torsor $\beta$ over $\bQ$ equipped with a trivialization $\beta$ over $\bA_f^p$ (resp. over $\bA_f$).
Then $G':=\Aut(\xi)$ is an inner form of $G$ (so the dual group of $G$ and $G'$ are canonically identified), equipped with an isomorphism $\theta: G(\bA^p_f)\cong G'(\bA^p_f)$ (resp. $\theta: G(\bA_f)\cong G'(\bA_f)$).
\end{dfn}

Now let $(G,X)$ and $(G',X')$ be two Shimura data, with $G'$ a prime-to-$p$ trivialized inner form of $G$.
Via $\theta$, one can transport $K^p\subset G(\bA_f^p)$ to an open compact subgroup ${K'}^p\subset G'(\bA_f^p)$. We identify the prime-to-$p$ (derived) Hecke algebra 
$H_{K^p,\La}$ (defined in the same way as in \eqref{E: derived hecke global}) with $H_{{K'}^p,\La}$ and simply write them as $H_{K^p,\La}$.
Let $K'_p\subset G'(\bQ_p)$ be an open compact subgroup and write $K'=K'_p{K'}^p$ for the corresponding level. 

We fix a quasi-split inner form $G^*_{\bQ_p}$ of $G_{\bQ_p}$ and $G'_{\bQ_p}$ equipped with a pinning $(B^*_{\bQ_p},T^*_{\bQ_p},e^*)$, and realize $G_{\bQ_p}$ as $J_b$ and $G'_{\bQ_p}$ as $J_{b'}$ for $b,b\in B(G^*_{\bQ_p})$. Under our assumption that $G$ and $G'$ are adjoint, such $b,b'$ exist and are unique. Then we have the conjectural coherent sheaf $\frakA_{K_p,\La}$ and $\frakA_{K'_p,\La}$ as in \eqref{E: sheaf A} on the stack $\Loc_{\CG,p}\otimes\La$ of local Langlands parameters for $G^*_{\bQ_p}$ over $\La$. 

\begin{conjecture}\label{C: geo JL}
For every choice of specialization map $\on{sp}: \overline{\eta}\to \overline{v}$, there is a natural map
\begin{multline}\label{E: geoJL}
R\Hom_{\Coh(\Loc_{\CG,p}\otimes\La)}\bigl(\widetilde{V_\mu}\otimes \frakA_{K_p,\La}, \widetilde{V_{\mu'}}\otimes \frakA_{K'_p,\La}\bigr)\\
\to R\Hom_{H_{K^p,\La}}\bigl(C_c(\Sh_{K}(G)_{\overline\eta}, \La[d_{\mu}]), C_c(\Sh_{K'}(G')_{\overline\eta},\La[d_{\mu'}])\bigr),
\end{multline}
compatible with compositions.  
In particular, there is an ($E_1$-)algebra homomorphism
\begin{equation}\label{E: S operator}
S:R\End_{\Coh(\Loc_{\CG,p}\otimes\La)}\bigl( \widetilde{V_\mu}\otimes \frakA_{K_p,\La}\bigr)\to R\End_{H_{K^p,\La}}\bigl(C_c(\Sh_{K}(G)_{\overline\eta},\La)\bigr),
\end{equation}
compatible with \eqref{E: geoJL}. In addition, the induced action
\begin{equation}\label{E: S action}
H_{K_p,\La}\stackrel{\eqref{E: Hecke alg isom}}{\cong} R\End(\frakA_{K_p,\La})\to R\End\bigl(\widetilde{V_\mu}\otimes \frakA_{K_p,\La}\bigr)\xrightarrow{S} R\End_{H_{K^p,\La}}\bigl(C_c(\Sh_{K}(G)_{\overline\eta},\La)\bigr)
\end{equation}
coincides with the natural Hecke action of $H_{K_p,\La}$ on $C_c(\Sh_{K}(G)_{\overline\eta},\La)$ (and therefore is independent of the specialization map $\on{sp}$).
\end{conjecture}
This conjecture would be a consequence of a Galois theoretic description of $C_c(\Sh_{K}(G)_{\overline\eta},\La)$ similar to Conjecture \ref{C: gen coh}, but its formulation does not require the existence of the stack of global Langlands parameters for $\bQ$. In any case, a step towards a Galois theoretical description of  $C_c(\Sh_{K}(G)_{\overline\eta},\La)$ might require Conjecture \ref{C: geo JL} as an input. We also remark that
as in the function field case, there is a more general version of such conjecture in \cite[Sect. 4.7]{Zhu-Conjecture}, allowing ``generalized level structures", so that the cohomology of Igusa varieties could appear.

The following theorem verifies the conjecture in special cases.
\begin{thm}
\label{introT: spectral action}
Suppose that the Shimura data $(G,X)$ and $(G',X')$ are of abelian type, with $G'$ a finitely trivialized inner form of $G$. Suppose that $G_{\bQ_p}$ is unramified (and therefore so is $G'_{\bQ_p}$).
\begin{enumerate}
\item\label{introT: spectral action-1} The map \eqref{E: geoJL} (and therefore \eqref{E: S operator}) exists when $\La=\bQ_\ell$ and  $K_p\subset G(\bQ_p)$ and $K'_p\subset G'(\bQ_p)$ are parahoric subgroups (in the sense of Bruhat-Tits).
\item\label{introT: spectral action-2} If $K_p$ is hyperspecial, then the map \eqref{E: geoJL} (and therefore \eqref{E: S operator}) exists when $\La=\bZ_\ell$, at least for underived $\Hom$ spaces. In addition, the
action of $H_{K_p}^{\cl}$ on $H_c^*(\Sh_{K}(G)_{\overline\eta},\La)$ via \eqref{E: S action} coincides with the natural action of $H_{K_p}^{\cl}$. 
\end{enumerate}
\end{thm}
Part \eqref{introT: spectral action-1} is proved in \cite{Xiao-Zhu-Cycle, Hemo-Zhu-Unipotent}. The proof contains two ingredients. One is the construction of physical correspondences between mod $p$ fibers of $\Sh_K(G)$ and $\Sh_{K'}(G')$ by \cite{Xiao-Zhu-Cycle} (this is where we currently need to assume that $G$ and $G'$ are unramified at $p$). The other ingredient is Theorem \ref{T: tame local Langlands} (and therefore requires $\La=\bQ_\ell$). When $K_p$ is hypersepcial, one can work with $\bZ_\ell$-coefficient, as (the underived version of) \eqref{E: local S-operator} exists for $\bZ_\ell$-coefficient thanks to \cite{Yu-Jacquet-Langlands-Higher-Weights}. In fact, in this case one can allow non-trivial local systems on the Shimura varieties (see  \cite{Yu-Jacquet-Langlands-Higher-Weights}).
The last statement is known as the $S=T$ for Shimura varieties. The case when $d_{\mu}=\dim\Sh_K(G)=0$ is contained in \cite{Xiao-Zhu-Cycle}. The general case is proved in \cite{Wu-S=T, Zhu-S=T} using foundational works from \cite{Scholze-Berkeley, Fargues-Scholze-Geometrization}. 

\section{Applications to arithmetic geometry}\label{Application}
Besides the previously mentioned directly applications of (ideas from) geometric Langlands to the classical Langlands program, we discuss some further arithmetic applications, mostly related to Shimura varieties and based on the author's works.
We shall mention that there are many other remarkable applications of (ideas of) geometric Langlands to arithmetic problems, such as \cite{Gaitsgory-deJong, HNY-Kloosterman, Yun-motive-exceptional-group, Xu-Zhu-Bessel, LLLM}, to name a few. 

\subsection{Local models of Shimura varieties}\label{S: local model}
The theory of integral models of Shimura varieties (with parahoric level) started (implicitly in the work of Kronecker) with understanding of the mod $p$ reduction of elliptic modular curves with $\Ga_0(p)$-level. We discuss a small fraction of this theory concerning \'{e}tale local structures of these integral models via the theory of local models. The recent developments of the theory of local models are greatly influenced by the geometric Langlands program.

We use notations from Sect. \ref{geometric JL} for Shimura varieties (but we do not assume that $G$ is of adjoint type in this subsection). Let $(G,X)$ be a Shimura datum and $K$ a chosen level with $K_p=\mG(\bZ_p)$ for some parahoric group scheme $\mG$ (in the sense of Bruhat-Tits) of $G_{\bQ_p}$ over $\bZ_p$. Then for a place $v$ of $E$ over $p$, a local model diagram is a correspondence of quasi-projective schemes over $\mO_{E_v}$
\begin{equation}\label{E: local model diagram}
\scrS_K\leftarrow\widetilde{\scrS}_K\xrightarrow{\tilde\varphi} M^{\loc}_{\mG},
\end{equation}
where $\scrS_K$ is an integral model of $\Sh_K(G)$ over $\mO_{E_v}$, $\widetilde{\scrS}_K$ is a $\mG_{\mO_{E_v}}$-torsor over $\scrS_K$, $M^{\loc}_{\mG}$ is the so-called local model, which is a flat projective scheme over $\mO_{E_v}$ equipped with a $\mG_{\mO_{E_v}}$-action, and $\tilde{\varphi}$ is a $\mG_{\mO_{E_v}}$-equivariant smooth morphism of relative dimension $\dim G$.
Therefore, $M_{\mG}^{\loc}$ models \'{e}tale local structure of $\scrS_K$. On the other hand, the existence of $\mG_{\mO_{E_v}}$-action on $M_{\mG}^\loc$ makes it easier than $\scrS_K$ to study. 

The original construction of local models is based on realization of a parahoric group scheme as (the neutral connected component of) the stabilizer group of a self-dual lattice chain in a vector space (over a division algebra over $F$) with a bilinear form, e.g. see  \cite{PRS-Survey} for a survey and references. This approach is somehow ad hoc and is limited the so-called (P)EL (local) Shimura data. A new approach, based on the construction of an $\bZ_p$-analogue of the stack $\Hk_{\mG,D}$ from Sect. \ref{Tamely local}, was systematically introduced in \cite{Pappas-Zhu-Local-Model} (under the tameness assumption of $G$ which was later lifted in  \cite{Levin-Weil-Restriction, Lourenco}). In \emph{loc. cit.} the construction of such $\bZ_p$-analogue (or rather the corresponding Beilinson-Drinfeld type affine Grassmannian $\Gr_{\mG,\bZ_p}$ over $\bZ_p$) is based on the construction of certain ``two dimensional parahoric" group scheme $\widetilde\mG$ over $\bZ_p[\varpi]$ whose restriction along $\bZ_p[\varpi]\xrightarrow{\varpi\mapsto p}\bZ_p$ recovers $\mG$. (See \cite{Zhu-ICCM} for a survey.)
A more direct construction of a different $p$-adic version of such affine Grassmannian $\Gr_{\mG,\Spd\bZ_p}$ was given in \cite{Scholze-Berkeley} in the analytic perfectoid world. In either case, the local model is defined as the flat closure of the Schubert variety  in the generic fiber corresponding to $\mu$. In addition, the recent work \cite{AGLR} shows that the two constructions agree.
The following theorem from  \cite{AGLR} is the most up-to-date result on the existence of local models and about their properties.
\begin{thm}
Let  $G$ be a connected reductive group over a $p$-adic field $F$. Except the odd unitary case when $p=2$ and triality case when $p=3$,
for every parahoric group scheme $\mG$ of $G$ over $\mO$, and a conjugacy class of minuscule cocharacters $\mu$ of $G$ defined over a finite extension $E/F$ of $F$,
there is a normal flat projective scheme $M_{\mG,\mu}^{\loc}$ over $\mO_E$, equipped with a $\mG_{\mO_E}$-action such that $M_{\mG,\mu}^{\loc}\otimes E$ is $G_E$-equivariantly isomorphic to the partial flag variety $\Fl_{\mu}$ of $G_E$ corresponding to $\mu$, and that $M_{\mG}^{\loc}\otimes k_E$ is $(\mG\otimes k_E)$-equivariantly isomorphic to the (canonical deperfection of the) union over the $\mu$-admissible set of Schubert varieties in $LG/L^+\mG\otimes k_E$. In addition, $M_{\mG}^{\loc}$ is
normal, Cohen-Macaulay and  each of its geometric irreducible components in its special fiber  is normal and Cohen-Macaulay. 
\end{thm}

We end this subsection with a few remarks.
\begin{rmk}
\begin{enumerate}
\item Once the local model diagram \eqref{E: local model diagram} is established, this theorem also gives the corresponding properties of the integral models of Shimura varieties. 

\item A key ingredient in the study of special fibers of local models is the coherence conjecture by Pappas-Rapoport \cite{Pappas-Rapoport-Ramified-Unitary}, proved in \cite{Zhu-Coherence-Conjecture} (and the proof uses the idea of fusion).

\item 
One important motivation/application of the theory of local models is the Haines-Kottwitz conjecture \cite{Haines-Test-Function}, which predicts certain central element in the parahoric Hecke algebra $H_{K_p}^{\cl}$ should be used as the test function in the trace formula computing the Hasse-Weil zeta function of $\Sh_{K}(G)$. As mentioned in Sect. \ref{Tamely local}, this conjecture motivated 
Gaitsgory's central sheaf construction \eqref{E: central sheaf}. With the local Hecke stack $\Hk_{\mG,\bZ_p}$ over $\bZ_p$ constructed (either the version from \cite{Pappas-Zhu-Local-Model} or from \cite{Scholze-Berkeley}), one can mimic the construction \eqref{E: central sheaf} in mixed characteristic to solve the Kottwitz conjecture. Again, see \cite{AGLR} for the up-to-date result.
\end{enumerate}
\end{rmk} 

\subsection{The congruence relation}\label{S: Congruence}
We use notations and (for simplicity) keep assumptions from Sect. \ref{geometric JL} regarding Shimura varieties.
Let $(G,X)$ be a Shimura datum abelian type, and let $K$ be a level such that $K_p$ is hyperspecial. Let $v\mid p$ be the place of $E$. Then $\Sh_K(G)$ has a canonical integral model $\scrS_K$ defined over $\mO_{E,(v)}$ (\cite{Kisin-Integral-Model}).
Let $\overline{\scrS}_K$ be its mod $p$ fiber, which is a smooth variety defined over the residue field $k_v$ of $v$. Let $\sigma_v$ denote the geometric Frobenius in $\Ga_{k_v}$. Theorem \ref{introT: spectral action} gives an action of $\End_{\Loc_{\CG,p}^{\ur}}(\widetilde{V_\mu})$ on $H^*_c(\overline\scrS_{K,\overline{k_v}},\bZ_\ell)$, which as we shall see has significant consequences. 

The congruence relation conjecture (also known as the Blasius-Rogawski conjecture), generalizing the classical Eichler-Shimura congruence relation $\Fr_p=T_p+V_p$ for modular curves,  predicts that in the Chow group of $\overline{\scrS}_K\times \overline{\scrS}_K$, the Frobenius endomorphism of $\overline{\scrS}_K$ satisfies a polynomial whose coefficients are mod $p$ reduction of certain Hecke correspondences. 
Theorem \ref{introT: spectral action}, together with \cite[Sect. 6.3]{Xiao-Zhu-Invariants}, implies this conjecture at the level of cohomology.

For every representation $V$ of ${}^c(G_{\bQ_p})$, its character $\chi_V$ (regarded as a $\hat{G}$-invariant function on $\CG|_{d=(p,\sigma_p)}$) gives an element $h_V\in H_{G(\bZ_p)}^{\cl}$ via the Satake isomorphism \eqref{Sat isom}.

\begin{thm}
The following identity 
\begin{equation}\label{E: congruence relation} 
\sum_{i=0}^{n}(-1)^j h_{\chi_{\wedge ^j V}} \sigma_v^{\dim V-j}=0
\end{equation}
holds in $\End (H^*_c(\overline{\scrS}_{K ,\overline{k_v}},\bZ_\ell))$, where $V=\on{Ind}_{{}^c(G_{E_v})}^{{}^c(G_{\bQ_p})}V_\mu$ is the tensor induction of $V_\mu$.
\end{thm}

Indeed, by \cite[Sect. 6.3]{Xiao-Zhu-Invariants}, such equality holds with $h_{\chi_{\wedge ^i V}}$ replaced by $S(\chi_{\wedge ^i V})$, where $S$ is from Theorem \ref{introT: spectral action} \eqref{introT: spectral action-1}. Then Part \eqref{introT: spectral action-2} of that theorem allows one to replace $S(\chi_{\wedge ^i V})$ by $h_{\chi_{\wedge ^i V}}$.  This approach to \eqref{E: congruence relation} is the Shimura variety analogue of V. Lafforgue's approach to the Eichler-Shimura relation for $\Sht_K(G)$ \cite{Lafforgue-Automorphic-to-Galois}. Traditionally, there is another approach to the congruence relation conjecture for Shimura varieties by directly studying reduction mod $p$ of Hecke operators, starting from \cite{Faltings-Chai-Abelian} for the Siegel modular variety case. See \cite{Lee-ES-Relation} for the latest progress and related references. This approach would give \eqref{E: congruence relation}  at the level of algebraic correspondences.

Now suppose $(G,X)=(\Res_{F^+/\bQ}(G_0)_{F^+}, \prod_{\varphi :F^+\to \bR} X_0)$, where $(G_0,X_0)$ is a Shimura datum and $F^+$ is a totally real field. As before, let $p$ be a prime such that $K_p$ is hyperspecial. In particular, $p$ is unramified in $F^+$. 
In addition, for simplicity we assume that $G_{0,\bQ_p}$ is split (so for a place $v$ of $E$ above $p$, $E_v=\bQ_p$). We let $\bF$ denote an algebraic closure of $\bF_p$.
Let $\{w_i\}_i$ be the set of primes of $F^+$ above $p$, and let $k_i$ denote the residue field of $w_i$. For each $i$, we also fix an embedding $\rho_i: k_i\to\bF$. 
Then there is a natural map 
\[
\prod_i(\bZ^{f_i}\rtimes \frakS_{f_i})\to \End_{\Loc_{\CG,p}^{\ur}}(\widetilde{V_\mu}),
\] 
where $\frakS_{f_i}$ is the permutation group on $f_i$ letters.
Together with Theorem \ref{introT: spectral action}, one obtains the following result (\cite{Xiao-Zhu-Cycle}).
\begin{thm}
There is an action of $\prod_i(\bZ^{f_i}\rtimes \frakS_{f_i})$ on $H_c^*(\overline\scrS_{K,\overline{\bF}},\bZ_\ell)$ such that action of $\sigma_p$ factors as $\sigma_p=\prod_i \sigma_{p,i}$, where $\sigma_{p,i}=((1,0,\ldots,0), (12\cdots f_i))\in \bZ^{f_i}\rtimes\frakS_{f_i}$.
Each $\sigma_{p,i}^{f_i}$ satisfies a polynomial equation similar to \eqref{E: congruence relation}.
\end{thm}
This theorem gives some shadow of the plectic cohomology conjecture of Nekov\'a\v{r}-Scholl \cite{Nekovar-Scholl-Plectic}.

\subsection{Generic Tate cycles on mod $p$ fibers of Shimura varieties}\label{S: Tate}
In \cite{Xiao-Zhu-Cycle}, 
we applied Theorem \ref{introT: spectral action} to verify ``generic" cases of Tate conjecture for the mod $p$ fibers of many Shimura varieties. We use notations and (for simplicity) keep assumptions from Sect. \ref{S: Congruence}. Let $(\overline{\mathscr S}_{K,\overline{k_v}})^\pf$ denote the perfection of $\overline{\mathscr S}_{K,\overline{k_v}}$ (i.e. regard it as a perfect presheaf over $\Aff_{\overline{k_v}}^{\pf}$), then by attaching to every point of $\overline{\mathscr S}_{K,\bar{k}}$ an $F$-isocrystal with $G$-structure (see \cite{Kisin-Integral-Model, Xiao-Zhu-Cycle}), one can define
the so-called Newton map 
\[
\on{Nt}: (\overline{\mathscr S}_{K,\overline{k_v}})^\pf\to \kot(G_{\bQ_p})_{\overline{k_v}}.
\] 
Then the Newton stratification of $\kot(G_{\bQ_p})_{\overline{k_v}}$ (see Sect. \ref{S: Local Langlands})
induces a stratification of $\overline{\mathscr S}_{K,\overline{k_v}}$ by locally closed subvarieties. It is known that the image of $\on{Nt}$ contains a unique basic element $b$ and the corresponding subvarieties in $\overline{\mathscr S}_{K,\overline{k_v}}$ is closed, called the basic Newton stratum, and denoted by $\overline{\scrS}_b$.

Let $m$ be the order of the action of the geometric Frobenius $\sigma_p$ on $\xch(\hat{T})$. Let
\[
\Lambda^\Tate_p=\Big\{\la\in \xch(\hat T)\; \Big|\; \sum_{i=0}^{m-1} \sigma_p^i(\la)=0\Big\}\subset \xch(\hat T).
\]
For a representation $V$ of $\hat G_{\bQ_\ell}$ and $\la\in\xch(\hat T)$, let $V(\la)$ denote the $\la$-weight subspace of $V$ (with respect to $\hat T$), and let
\[
V^\Tate=\bigoplus_{\la\in\Lambda^\Tate_p} V(\la).
\]

We are in particular interested in the condition $V^\Tate_\mu\neq 0$.
As explained in the introduction of \cite{Xiao-Zhu-Cycle}, under the conjectural Galois theoretic description of the cohomology of the Shimura varieties (analogous to Conjecture \ref{C: gen coh}), for a Hecke module $\pi_f$ whose Satake parameter at $p$ is general enough, certain multiple $a(\pi_f)$ of
the dimension of this vector space should be equal to the dimension of the space of Tate classes in the $\pi_f$-component of the middle dimensional compactly-supported cohomology of $\overline{\scrS}_{K,\overline{k_v}}$. 
In addition, this space is usually large. For example, in the case $G$ is an odd (projective) unitary  group of signature $(i,n-i)$ over a quadratic imaginary field, the dimension of this space at an inert prime is $\begin{pmatrix}\frac{n+1}{2}\\ i\end{pmatrix}$.

For a (not necessarily irreducible) algebraic variety $Z$ of dimension $d$ over an algebraically closed field, let $H^{\on{BM}}_{2d}(Z)(-d)$ denote the $(-d)$-Tate twist of the top degree Borel--Moore homology, which is the vector space spanned by the irreducible components of $Z$. Now let $X$ be
a smooth variety of dimension $d+r$ defined over a finite field $k$ of $q$ elements, and let $Z \subseteq X_{\overline k}$ be a (not necessarily irreducible) projective subvariety of dimension $d$. There is the cycle class map
\[
\cl: H^{\on{BM}}_{2d}(Z)(-d)\to \bigcup_{j\geq 1}H_{c}^{2d}(X_{\overline k}, \bQ_\ell(d))^{\sigma_q^j}=:T^d_\ell(X).
\]
\begin{thm}
\label{T:main theorem}
We write $d_{\mu}=\dim X=2d$ and $r=\dim V_\mu^{\Tate}$.
\begin{enumerate}
\item\label{T:main theorem-1} The basic Newton stratum $\overline{\mathscr S}_b$ of $\overline{\scrS}_{K,\overline{k_v}}$ is pure of dimension $d$. In particular, $d$ is always an integer. In addition,
there is an $H_{K,\bQ_\ell}$-equivariant isomorphism
\[
H^{\on{BM}}_{2d}(\overline{\mathscr S}_{b})(-d)\cong C(G'(\bQ)\backslash G'(\bA_f)/K,\bQ_\ell)^{\oplus r},
\]
where $G'$ is the finitely trivialized inner form of $G$ with $G'_\bR$ is compact.
\item\label{T:main theorem-2}
Let $\pi_f$ be an irreducible module of $H_{K,\Ql}$, and let 
\[
H^{\on{BM}}_{2d}(\overline{\mathscr S}_{b})[\pi_f]=\Hom_{H_{K,\Ql}}(\pi_f,H^{\on{BM}}_{2d}(\overline{\mathscr S}_{b})(-d)_{\Ql})\otimes \pi_f
\] 
be the $\pi_f$-isotypical component. Then the cycle class map
\[
\cl:H^{\on{BM}}_{2d}(\overline{\mathscr S}_{b})(-d)\to T^d_\ell(\overline{\mathscr S}_{K})
\]
restricted to $H^{\on{BM}}_{2d}(\overline{\mathscr S}_{b})[\pi_f]$ is injective if the Satake parameter of $\pi_{f,p}$ (the component of $\pi_f$ at $p$) is $V_\mu$-general.
\item\label{T:main theorem-3} 
Assume that $\Sh_K(G)$ is (essentially) a quaternionic Shimura variety or a Kottwitz arithmetic variety. Then the $\pi_f$-isotypical component of the cycle class map is surjective to $T^d_\ell(\overline{\scrS_K})[\pi_f]$ if the Satake parameter of $\pi_{f,p}$ is strongly $V_\mu$-general. In particular, the Tate conjecture holds for these $\pi_f$.
\end{enumerate}
\end{thm}
\begin{rmk}
\begin{enumerate}
\item For a representation $V$ of $\hat G$, the definitions of ``$V$-general" and ``strongly $V$-general" Satake parameters can be found in \cite[Definition 1.4.2]{Xiao-Zhu-Cycle}. Regular semisimple elements in $\CG|_{d=(p,\sigma_p)}$ are always $V$-general, but not the converse. See \cite[Remark 1.4.3]{Xiao-Zhu-Cycle}.
\item Some special cases of the theorem were originally proved in \cite{Helm-Tian-Xiao,Tian-Xiao-Tate}.
\end{enumerate}
\end{rmk}
The proof of this theorem relies on several different ingredients. Via the Rapoport-Zink uniformization of the basic locus of a Shimura variety, Part \eqref{T:main theorem-1} can be reduced a question about irreducible components of certain affine Deligne-Lusztig varieties, which was studied in \cite[\S 3]{Xiao-Zhu-Cycle}.
The most difficult is Part \eqref{T:main theorem-2}, which we proved by calculating the intersection numbers among all $d$-dimensional cycles in $\overline\scrS_b$. These numbers  can be encoded in an $r\times r$-matrix with entries in $H^{\cl}_{K_p}$. In general, it seems hopeless to calculate this matrix directly and explicitly. 
However, this matrix can be understood as the composition of certain morphisms in $\Coh(\Loc_{\CG,p}^{\ur})$. Namely, first we realize $G'(\bQ)\backslash G'(\bA)/K$ as a Shimura set with $\mu'=0$ its Shimura cocharacter.
Then using Theorem \ref{introT: spectral action} (and the Satake isomorphism \eqref{E: derived Satake}), this matrix can be calculated as
\[ 
\Hom_{\Coh(\Loc_{\CG,p}^{\ur})}(\mO,\widetilde{V_{\mu}})\otimes  \Hom_{\Coh(\Loc_{\CG,p}^{\ur})}(\widetilde{V_{\mu}},\mO)\to  \Hom_{\Coh(\Loc_{\CG,p}^{\ur})}(\mO,\mO)\cong H^{\cl}_{K_p}\otimes\bQ_\ell.
\]
Then one needs to determine when this pairing is non-degenerate, which itself is an interesting question in representation theory, whose solution relies on the study of the Chevellay's restriction map for vector-valued functions. The determinant of this matrix was calculated in \cite{Xiao-Zhu-Invariants}. Finally, Part \eqref{T:main theorem-3} was proved by comparing two trace formulas, the Lefschetz trace formula for $G$ and the Arthur-Selberg trace formula for $G'$. 

\begin{ex}\label{Tate cycle on U(1,2n)}
Let $G=\on{U}(1,2r)$ be the unitary group\footnote{This is not an adjoint group so the example is not consistent with our assumption. But it is more convenient for the discussion here. The computations are essentially the same.} of $(2r+1)$-variables associated to an imaginary quadratic extension $E/\bQ$, whose signature is $(1,2r)$ at infinity. It is equipped with a standard Shimura datum, giving a Shimura variety (after fixing a level $K\subset G(\bA_f)$). In particular if $r=1$, this is (essentially) the Picard modular surface. Let $p$ be a prime inert in $E$ such that $K_p$ is hyperspecial.  In this case $\overline{\scrS}_b$ is a union of certain Deligne-Lusztig varieties, parametrized by $G'(\bQ)\backslash G'(\bA_f)/K$, where $G'=\on{U}(0,2r+1)$ that is isomorphic to $G$ at all finite places. The intersection pattern of these cycles inside  $\overline{\scrS}_b$ were (essentially) given in \cite{Vollaard-Wedhorn-GU(n1)} but the intersection numbers between these cycles are much harder to compute. In fact we do not know how to compute them directly for general $r$, except applying Theorem \ref{introT: spectral action} to this case. (The case $r=1$ can be handled directly.)

We have $\hat G=\GL_{2r+1}$ on which $\sigma_p$ acts as $A\mapsto J(A^T)^{-1}J$, where $J$ is the anti-diagonal matrix with all entries along the anti-diagonal being $1$. The representation $V_\mu$ is the standard representation of $\GL_{2r+1}$. One checks that $\dim V_\mu^{\Tate}=1$ (which is consistent with the above mentioned parameterization of irreducible components of $\overline{\scrS}_b$ by $G'(\bQ)\backslash G'(\bA_f)/K$). 
We identify the weight lattice of $\hat G$ as $\bZ^{2r+1}$ as usual. Then $\Hom_{\Coh(\Loc_{\CG,p}^{\ur})}(\mO,\widetilde{V_{\mu}})$ is a free rank one module over $\Hom_{\Coh(\Loc_{\CG,p}^{\ur})}(\mO,\mO)=H_{K_p}^{\cl}\otimes\bQ_\ell$. Then a generator $\mathbf a_{\on{in}}$ induces an $H_{K,\bQ_\ell}$-equivariant homomorphism
$$S(\mathbf a_{\on{in}}): C(G'(\bQ)\backslash G'(\bA_f)/K)\to H_c^{2r}(\overline{\scrS}_{K,\overline{k_v}},\bQ_\ell(r)),$$  
realizing the cycle class map of $\overline\scrS_b$ (up to a multiple). The module  $\Hom_{\Coh(\Loc_{\CG,p}^{\ur})}(\widetilde{V_{\mu}},\mO)$ is also free of rank one over $H_{K_p,\bQ_\ell}$. For a chosen generator $\mathbf a_{\on{out}}$, the composition
\[S(\mathbf a_{\on{out}})\circ S(\mathbf a_{\on{in}})=S(\mathbf a_{\on{out}}\circ\mathbf a_{\on{in}})\]
calculates the intersection matrix of those cycles from the irreducible components of $\overline\scrS_b$. 

The element $h:=\mathbf a_{\on{out}}\circ\mathbf a_{\on{in}} \in H_{K_p,\bQ_\ell}$ was explicitly computed in \cite[Example 6.4.2]{Xiao-Zhu-Invariants} (up to obvious modification and also via the Satake isomorphism \eqref{Sat isom classical}). Namely,
\begin{equation}
\label{E:Hecke operator in Tpj}
h=p^{r(r+1)}\sum_{i=0}^{r} (-1)^i(2i+1)p^{(i-r)(r+i+1)}\sum_{j=0}^{r-i} \begin{bmatrix} 2r+1-2j \\ r-i- j\end{bmatrix}_{t=-p}T_{p,j}.
\end{equation}
Here, $T_{p,j}=1_{K_p\la_j(p)K_p}$, with $\la_i=(1^i,0^{2r-2i+1},(-1)^i)$, and $\begin{bmatrix} n \\ m\end{bmatrix}_t$ is the
 $t$-analogue of the binomial coefficient given by
\[
[0]_t=1,\quad [n]_t=\frac{t^n-1}{t-1},\quad [n]_t!=[n]_t[n-1]_t\cdots[1]_t,\quad \begin{bmatrix} n \\ m\end{bmatrix}_t=\frac{[n]_t!}{[n-m]_t![m]_t!}.
\]
In other words, the intersection matrix of cycles in $\overline\scrS_b$ in this case is calculated by the Hecke operator \eqref{E:Hecke operator in Tpj}. 

On interesting consequence is this computation is the following consequence on the intersection theory of the finite Deligne-Lusztig varieties, for which we do not know a direct proof.
Let $W$ be a $(2r+1)$-dimensional non-degenerate hermitian space over $\bF_{p^2}$.
Consider the following $r$-dimensional Deligne-Lusztig variety
\[
\mathrm{DL}_{r}: = \big\{ H \subset W \textrm{ of dimension }r\; |\; H \subseteq (H^{(p)})^\perp \big\},
\]
where $H^{(p)}$ the pullback of $H$ along the Frobenius.
Let $\mH$ denote the corresponding universal subbundle of rank $r$. Let $\mE = \mH^{(p)} \otimes \big( (\mH^{(p)})^\perp / \mH \big)$. Then we have
\begin{equation}\label{E: chern class on DL}
\int_{\mathrm{DL}_r} c_r( \mE)  = \sum_{i = 0}^r (-1)^i(2i+1) p^{i^2+i} \begin{bmatrix} 2r+1 \\ r-i\end{bmatrix}_{t=-p}.
\end{equation}
\end{ex}

\subsection{The Beilinson-Bloch-Kato conjecture for Rankin-Selberg motives}\label{S: BBK}
Let $M$ be a rational pure Chow motive of weight $-1$ over a number field $F$. The Beilinson-Bloch-Kato conjecture, which is a far reaching generalization of the Birch and Swinnerton-Dyer conjecture, predicts deep relations between certain  algebraic, analytic, and cohomological invariants attached to $M$:
\begin{itemize}
\item the rational Chow group $\CH(M)^0$ of homologically trivial cycles of $M$;
\item the $L$-function  $L(s,M)$ of $M$;
\item the Bloch-Kato Selmer group $H^1_f(F, H_{\ell}(M))$ of the $\ell$-adic realization $H_{\ell}(M)$ of $M$. 
\end{itemize}
The Beilinson-Bloch conjecture predicts an equality 
\[
\dim_{\bQ}\CH(M)^0 = \on{ord}_{s=0}L(s,M)
\] 
between the dimension of $\CH(M)^0$ and the vanishing order of the $L$-function at the central point, while the Bloch-Kato conjecture predicts 
\[
\on{ord}_{s=0}L(s,M)= \dim_{\bQ_\ell} H^1_f(F, H_{\ell}(M)).
\] 
In addition, the so-called $\ell$-adic Abel-Jacobi map
\[
\on{AJ}_\ell: \on{CH}(M)^0\otimes \bQ_\ell\to H^1_f(F, H_{\ell}(M))
\] 
should be an isomorphism.

This conjecture seems to be completely out of reach at the moment. E.g. for a general motive it is still widely open whether the $L$-function has a meromorphic continuation to the whole complex plane so that the vanishing order of $L(s,M)$ at $s=0$ makes sense. (This would follow from the Galois-to-automorphic direction of the Langlands correspondence for number fields.) Despite of this, there have been many works testing this conjecture in various special cases, mostly for motives $M$ of small rank.
In the work \cite{LTXZZ-BBK}, we verify certain cases of the above conjecture for Rankin-Selberg motives, which consist of a sequence of motives of arbitrarily large rank. 

We assume that $F/F^+$ is a (non-trivial) CM extension with $F^+$ totally real in the sequel. 
\begin{thm}\label{RS-Sym-Elliptic}
Let $A_1,A_2$ be two elliptic curves over $F^+$. Assume that
\begin{enumerate}
\item\label{RS-Sym-Elliptic-1} $\End_{\overline F}A_i=\bZ$ and $\Hom_{\overline F}(A_1,A_2)=0$;
\item\label{RS-Sym-Elliptic-2} $\Sym^{n-1}A_1$ and $\Sym^nA_2$ are modular;
\item\label{RS-Sym-Elliptic-3} $F^+\neq \bQ$ if $n\geq 3$.
\end{enumerate}
Under these assumption, if $L(n, \Sym^{n-1}A_1\times \Sym^nA_2)\neq 0$, then for almost $\ell$, 
\[
\dim_{\bQ_\ell} H^1_f(F,  \Sym^{n-1}V_\ell(A_1)\otimes\Sym^nV_\ell(A_2)(1-n))=0.
\]
Here $V_\ell(A_i)$ denotes the rational Tate module of $A_i$ as usual.
\end{thm}

This theorem is in fact a consequence of a more general result concerning Bloch-Kato Selmer groups of Galois representations associated to certain Rankin-Selberg automorphic representations, which we now discuss.

Recall that for an irreducible regular algebraic conjugate self-dual cuspidal  (RACSDC) automorphic representation $\Pi$ of $\GL_{n}(\bA_F)$, one can attach a compatible system of Galois representations $\rho_{\Pi,\la}:\Ga_F\to \GL_n(E_\la)$, where $E\subset \bC$ is a large enough number field and $\la$ is a prime of $E$ (see \cite{Chenevier-Harris-Galois-Rep}). Such $E$ is called a strong coefficient field of $\Pi$, which in the situation considered below can be taken as the number field generated by Hecke eigenvalues of $\Pi$.

\begin{thm}\label{BBK for RS}
Suppose that $F^+\neq\bQ$ if $n\geq 3$.
Let $\Pi_n$ (resp. $\Pi_{n+1}$) be an RACSDC automorphic representation of $\GL_{n}(\bA_F)$ (resp. $\GL_{n+1}(\bA_F)$) with trivial infinitesimal character. Let $E\subseteq\bC$ be a strong coefficient field for both $\Pi_n$ and $\Pi_{n+1}$. Let $\lambda$ be an admissible prime of $E$ with respect to $\Pi:=\Pi_0\times \Pi_1$. Let $\rho_{\Pi,\la}:=\rho_{\Pi_n,\la}\otimes_{E_\la} \rho_{\Pi_{n+1},\la}$.
\begin{enumerate}
\item\label{BBK for RS-1} If the Rankin-Selberg $L$-value $L(\frac{1}{2},\Pi)\neq 0$\footnote{Here we use the automorphic normalization of the $L$-function.}, then $H^1_f(F,\rho_{\Pi,\lambda}(n))=0$. 
\item\label{BBK for RS-2} If certain element $\Delta_\la\in H^1_f(F,\rho_{\Pi,\lambda}(n))$ (to be explained below) is non-zero, then $H^1_f(F,\rho_{\Pi,\lambda}(n))$ is generated by $\Delta_\la$ as an $E_\la$-vector space.
\end{enumerate}
\end{thm}

The notion of admissible primes appearing in the theorem consists of a long list of assumptions, some of which are rather technical. Essentially, it guarantees that the Galois representation $\rho_{\Pi,\lambda}$ has a well defined $\mO_{E,\la}$-lattice (still denoted by $\rho_{\Pi,\lambda}$ in the sequel) and the reduction mod $\la$ representation is suitably large and contains certain particular elements. (This is also related to the notion of $V$-general from Theorem \ref{T:main theorem}.)
Fortunately, in some favorable situations, one can show that all but finitely many primes are admissible. 
For example, this is the case considered in Theorem \ref{RS-Sym-Elliptic}. For another case in pure automorphic setting, see \cite[Thm. 1.1.7]{LTXZZ-BBK}.

The proof of the theorem uses several different ingredients. The initial step for Case \eqref{BBK for RS-1} is to translate the analytic condition $L(\frac{1}{2},\Pi)\neq 0$ into a more algebraic condition via the global Gan-Gross-Prasad (GGP) conjecture. Namely, the GGP conjecture predicts that in this case, there exist a pair of hermitian spaces $(V_n,V_{n+1})$ over $F$  that are totally positive definite at $\infty$, where $V_{n+1}=V_n\oplus Fv$ with $(v,v)=1$, and a
tempered cuspidal automorphic representation $\pi=\pi_n\times \pi_{n+1}$ of the product of unitary groups $G=U(V_n)\times U(V_{n+1})$, such that the period integral 
\[
[\Delta_H]: C^*_c(\Sh(G),E)[\pi]\to E
\]
does not vanish, where $H:=U(V_n)$ embeds into $G$ diagonally, which induces an embedding $\Delta_H: \Sh(H)\hookrightarrow \Sh(G)$ of corresponding Shimura varieties (in fact Shimura sets) with appropriately chosen level structures (here and later we omit level structures from the notations). We denote by $[\Delta_H]$ the homology class of $\Sh(G)$ given by $\Sh(H)$ and write $C^*_c(\Sh(G),E)[\pi]$ for the $\pi$-isotypical component of the cohomology (i.e. functions) of $\Sh(G)$. 
This conjecture was first proved in \cite{Zhang-GGP}
under some local restrictions which are too restrictive for arithmetic applications. Those restrictions are all lifted in our recent work through some new techniques in the study of trace formulae (\cite{BLZZ-GGP}). 

The strategy then is to construct, for every $m\geq 1$, (infinitely many) cohomology classes $\{\Theta_m^p\}_p\subset H^1(F, (\rho_{\Pi,\la}/\la^m)^*(1))$, where $p$ are appropriately chosen primes and $(-)^*(1)$ denotes the usual Pontryagin duality twisted by the cyclotomic character, such that the local Tate pairing at $p$ between $\Theta_m^p$ and Selmer classes of the Galois representation $\rho_{\Pi,\la}/\la^m$  is related to the above period integral. Then one can use Kolyvagin type argument (amplified in \cite{Liu-HZ-Cycle, LTXZZ-BBK}), with  $\{\Theta_m^p\}$ served as annihilators of the Selmer groups, to conclude.

The construction of $\Theta_m^p$ uses the diagonal embedding of Shimura varieties 
\[
\Delta_{H'}: \Sh(H')\hookrightarrow \Sh(G')
\]
where $H'\hookrightarrow G'$ are prime-to-$p$ trivialized (extended pure) inner forms of $H\subset G$ (see Definition \ref{D: prime-to-p trivial}). These Shimura varieties have parahoric level structures at $p$, and using the theory of local models (Sect. \ref{S: local model}) one can show that their integral models are poly-semistable at $p$ and compute the sheaf of nearby cycles on their mod $p$ fibers. 
Using many ingredients, including the understanding of (integral) cohomology of $\Sh(G')$ over $\overline{F}$, the computations from Example \ref{Tate cycle on U(1,2n)} (in particular \eqref{E:Hecke operator in Tpj} and \eqref{E: chern class on DL}), and the  Taylor-Wiles patching method \cite{LTXZZ-Rigidity}, one shows that $(\rho_{\Pi,\la}/\la^m)^*(1)$ does appear in the cohomology of $\Sh(G')$ (the so-called arithmetic level raising for $\Pi$), and that the diagonal cycle $\Delta_{H'}$, when localized at $(\rho_{\Pi,\la}/\la^m)^*(1)$, does give the desired class $\Theta_m^p$. We shall mention that this is consistent with conjectures in \S~\ref{S: Local Langlands} and \S~\ref{geometric JL}, as coherent sheaves on $\Loc_{\CG,p}\otimes \mO_E/\la^m$ corresponding to $c\on{-ind}_{K_p}^{G(\bQ_p)}(\mO_E/\la^m)$ and $c\on{-ind}_{K'_p}^{G'(\bQ_p)}(\mO_E/\la^m)$ are expected to be related exactly in this way. 

We could also explain the class $\Delta_\la$ appearing in Case \eqref{BBK for RS-2}. Namely, in this case we start with an embedding of Shimura varieties $\Delta_H: \Sh(H)\hookrightarrow\Sh(G)$, where $G$ is a product of unitary groups such that $\Pi$ descends to a tempered cuspidal automorphic representation $\pi$ appearing in the middle dimensional cohomology of $\Sh_G$. Then the $\pi$-isotypical component of the cycle $\Delta_H$ is homologous to zero and we let $\Delta_\la=\on{AJ}_\la(\Delta_H[\pi])$. The strategy to prove Case \eqref{BBK for RS-2} then is to reduce it to Case \eqref{BBK for RS-1} via some similar arguments as before.

%------

%------
% Insert acknowledgments and information
% regarding funding at the end of the last
% section, i.e., right before the bibliography.
%------

\quash{
\begin{ack}
We thank X.
\end{ack}

\begin{funding}
This work was partially supported by~\ldots
\end{funding}
}
%------
% Insert the bibliography.
%------


\begin{thebibliography}{99}

%------ Example for a paper in journal:
% \bibitem{article1}
% A.~Petrunin, Parallel transportation for Alexandrov space with curvature bounded below.
% \emph{Geom. Funct. Anal.} \textbf{8} (1998), no.~1, 123--148.

%------ Example for a book:
% \bibitem{book1}
% W.~P. Ziemer, \emph{Weakly differentiable functions}.
% Grad. Texts in Math. 120,  Springer, New York, 1989.

%------ Example for a paper in a book:
% \bibitem{incollection1}
% J.~S. Milne, Introduction to Shimura varieties.
% In \emph{Harmonic analysis, the trace formula, and Shimura varieties},
% edited by M.~W. Marcellin and E.~Giorgi, pp. 265--378,
% Clay Math. Proc. 4, Amer. Math. Soc., Providence, RI, 2005.

%------ Example for a preprint on arXiv:
% \bibitem{preprint1}
% D.~V. Nguyen, S.~K. Chilappagari, M.~W. Marcellin, and B.~Vasic,
% LDPC codes from latin squares free of small trapping sets,
% 2010, \href{http://arxiv.org/abs/1008.4177}{arXiv:1008.4177}.

%------ Example for a report:
% \bibitem{report1}
% J.~Schöberl, Commuting quasi-interpolation operators.
% Technical report isc-01-10-math, Texas A\&M University, 2001,
% \url{www.isc.tamu.edu/publications-reports/tr/0110.pdf}.

%------ Example for a thesis:
% \bibitem{thesis1}
% E.~Giorgi, \emph{The geometric universe}.
% Ph.D. thesis, University of Maryland, College Park, 2002. 

\bibitem{AGLR}J.~Ansch\"utz and I.~Gleason and J.~Louren\c{c}o and T.~Richarz, On the $p$-adic theory of local models. ArXiv: 2201.01234.

\bibitem{Arinkin-Gaitsgory-Singular-Support}D.~Arinkin and D.~Gaitsgory, Singular support of coherent sheaves, and the geometric Langlands conjecture. \emph{Selecta Math.} \textbf{21} (2015), 1--199.

\bibitem{AGKRRV}D.~Arinkin and D.~Gaitsgory and D.~Kazhdan and S.~Raskin and N.~Rozenblyum and Y.~Varshavsky, The stack of local systems with restricted variation and geometric Langlands theory with nilpotent singular support. ArXiv:2010.01906.

\bibitem{Beauville-Laszlo}A.~Beauville and Y.~Laszlo, Conformal blocks and generalized theta functions. \emph{Comm. Math. Phys.} \textbf{164} (1994), no.~2, 385--419. 

\bibitem{Beilinson-Drinfeld-Hitchin}A.~Beilinson and V.~Drinfeld, Quantization of Hitchin's integrable system and Hecke eigensheaves. \href{www.math.uchicago.edu/~mitya/langlands}.

\bibitem{BCHN}D.~Ben-Zvi and H.~Chen and D.~Helm and D.~Nadler, Coherent Springer theory and the categorical Deligne-Langlands correspondence. ArXiv:2010.02321.

\bibitem{Ben-Zvi-Nadler-Betti-Geometric-Langlands}D.~Ben-Zvi and D.~Nadler, Betti geometric Langlands. In \emph{Algebraic Geometry: Salt Lake City 2015}, edited by T.~de Fernex et al., pp. 3--41, Proc. Symp. Pure Math. 97.2, Amer. Math. Soc., Providence, RI, 2018.

\bibitem{BLZZ-GGP} R.~Beuzart-Plessis and Y.~Liu and W.~Zhang and X.~Zhu. Isolation of cuspidal spectrum, with application to the Gan-Gross-Prasad conjecture. \emph{Ann. Math.} \textbf{194} (2021), no.~2, 519--584.

\bibitem{Bezrukavnikov-Affine-Hecke}R.~Bezrukavnikov, On two geometric realizations of an affine Hecke algebra. \emph{Publ. Math. IH\'ES} \textbf{123}, (2016) 1--67.

\bibitem{Bezrukavnikov-Braverman-Geometric-Langlands-char-p}R.~Bezrukavnikov and A.~Braverman, Geometric Langlands
conjecture in characteristic p: The $GL_n$ case. \emph{Pure Appl. Math. Q.} \textbf{3} (2007), no.~1, Special Issue: In honor of Robert D. MacPherson, Part 3, 153--179.

\bibitem{Bezrukavnikov-Finkelberg-Derived-Satake}R.~Bezrukavnikov and M.~Finkelberg, Equivariant Satake category and Kostant-Whittaker reduction. \emph{Mosc. Math. J.} \textbf{8} (2008), no.~1, 39--72.

\bibitem{Bezrukavnikov-Travkin-Quantization}R.~Bezrukavnikov and R.~Travkin and T.-H.~Chen and X.~Zhu, Quantization of Hitchin integrable system via positive characteristic. ArXiv:1603.01327.

\bibitem{Bhatt-Scholze-Witt}B.~Bhatt and P.~Scholze, Projectivity of Witt vector affine Grassmannians. \emph{Invent. Math.} \textbf{209} (2017), no.~2, 329--423.

\bibitem{Chen-Zhu-Nonabelian-Hodge}T.-H.~Chen and X.~Zhu, Non-abelian Hodge theory for curves in characteristic $p$. \emph{Geometric and Functional Analysis} \textbf{25} (2015), 1706--1733.

\bibitem{Chen-Zhu-Geometric-Langlands-char-p}T.-H.~Chen and X.~Zhu, Geometric Langlands in prime characteristic. \emph{Comp. Math.} \textbf{153} (2017), 395--452.

\bibitem{Chenevier-Harris-Galois-Rep}G.~Chenevier and M.~Harris, Construction of automorphic Galois representations, II, \emph{Camb. J. Math.} \textbf{1} (2013), no.~1, 53--73.

\bibitem{DHKM}J.F.~Dat and D.~Helm and R.~Kurinczuk and G.~Moss, Moduli of Langlands parameters. ArXiv: 2009.06708.

\bibitem{DLYZ-Endoscopy}G.~Dhillon and Y.~Li and Z~Yun and X.~Zhu, in preparation. 

\bibitem{Donagi-Pantev-Hitchin}R.~Donagi and T.~Pantev, Langlands duality for Hitchin systems. \emph{Invent. Math.} \textbf{189} (2012), 653--735.

\bibitem{Drinfeld-GL2}V. G.~Drinfeld, Two-dimensional $\ell$-adic representations of the fundamental group of a curve over a finite field and automorphic forms on $\GL(2)$, \emph{Amer. J. Math.} \textbf{105} (1983), 85--114.

\bibitem{Emerton-Gee-Zhu}M.~Emerton and T.~Gee and X.~Zhu, in preparation.

\bibitem{Emerton-Zhu}M.~Emerton and X.~Zhu, in preparation.

\bibitem{Faltings-Loop}G.~Faltings, Algebraic loop groups and moduli spaces of bundles. \emph{J. Eur. Math. Soc.} \textbf{5} (2003), no.~1, 41--68.

\bibitem{Faltings-Chai-Abelian}G.~Faltings and C.L.~Chai, \emph{Degeneration of abelian varieties}. volume 22 of Results in Mathematics and Related Areas (3), Springer-Verlag, Berlin, 1990.

\bibitem{Fargues-Scholze-Geometrization}L.~Fargues and P.~Scholze, Geometrization of the local Langlands correspondence. ArXiv:2102.13459.

\bibitem{FGV-Geometric-Langlands-GLn}E.~Frenkel and D.~Gaitsgory and K.~Vilonen, On the geometric Langlands conjecture. \emph{J. Amer. Math. Soc.} \textbf{15} (2001), 367--417.

\bibitem{Gaitsgory-Central-Sheaf}D. Gaitsgory, Construction of central elements in the affine Hecke algebra via nearby cycles. \emph{Invent. math.} \textbf{144} (2001), 253--280.

\bibitem{Gaitsgory-deJong}D. Gaitsgory, On De Jong's conjecture. \emph{Isr. J. Math.} \textbf{157} (2007), 155--191. 

\bibitem{Haines-Test-Function}T.~Haines, The stable Bernstein center and test functions for Shimura varieties. In \emph{Automorphic forms and Galois representations. Vol. 2}, edited by
F.~Diamond, P.L.~Kassaei and M.~Kim, pp.118--186, London Math. Soc. Lecture Note Ser. 415, Cambridge Univ. Press, Cambridge, 2014.

\bibitem{Harris-Taylor-GLn}M.~Harris and R.~Taylor, \emph{The Geometry and Cohomology of Some Simple Shimura Varieties}.  Ann. Math. Studies 151, Princeton University Press, 2001.

\bibitem{HNY-Kloosterman}J.~Heinloth and B.C.~Ng\^o and Z.~Yun, Kloosterman sheaves for reductive groups. \emph{Ann. Math.} \textbf{177} (2013), 241--310.

\bibitem{Hellmann-Iwahori}E. Hellmann, On the derived category of the Iwahori-Hecke algebra. ArXiv:2006.03013.

\bibitem{Helm-Tian-Xiao}D.~Helm and Y.~Tian and L.~Xiao, Tate cycles on some unitary Shimura varieties mod $p$. \emph{Algebraic Number Theory} \textbf{11} (2017), 2213--2288.

\bibitem{HRS-Kunneth}T.~Hemo and T.~Richarz and J.~Scholbach, Constructible sheaves on schemes and a categorical K\"unneth formula. ArXiv:2012.02853.

\bibitem{Hemo-Zhu-Unipotent}T.~Hemo and X.~Zhu, Unipotent categorical local Langlands correspondence. Preprint.

\bibitem{Kazhdan-Lusztig-Iwahori}D.~Kazhdan and G.~Lusztig, Proof of the Deligne-Langlands conjecture for Hecke algebras.\emph{Invent. Math.} \textbf{87} (1987), 153--215.

\bibitem{Kisin-Integral-Model}M.~Kisin, Integral models for Shimura varieties of abelian type. \emph{J. Amer. Math. Soc.} \textbf{23} (2010), 967--1012.

\bibitem{Lafforgue-GLn}L.~Lafforgue, Chtoucas de Drinfeld et correspondance de Langlands. \emph{Invent. math.} \textbf{147} (2002), 1--241.

\bibitem{Lafforgue-Automorphic-to-Galois}V.~Lafforgue, Chtoucas pour les groupes r\'eductifs et param\`etrisation de Langlands globale. \emph{J. Amer. Math. Soc.} \textbf{31} (2018), 719--891.

\bibitem{Lafforgue-Zhu-Elliptic}V.~Lafforgue and X.~Zhu, D\'ecomposition au-dessus des param\`etres de Langlands elliptiques. ArXiv: 1811.07976.

\bibitem{Langlands-to-Weil}R.~Langlands, Letter to Andr\'e Weil. \href{publications.ias.edu/rpl/paper/43}.

\bibitem{Laumon-Geometric-Langlands}G.~Laumon, Correspondance de Langlands g\'eom\'etrique pour les corps de fonctions.\emph{Duke Math. J.} \textbf{54} (1987), 309--359.

\bibitem{LRS-GLn}G.~Laumon and M.~Rapoport and U.~Stuhler, $D$-elliptic sheaves and the Langlands correspondence. \emph{Invent. Math.} \textbf{113} (1993), 217--338.

\bibitem{LLLM}D.~Le and B.V.~Le Hung and B.~Levin and S. Morra, Local models for Galois deformation rings and applications. ArXiv:2007.05398.

\bibitem{Lee-ES-Relation}S.Y.~Lee, Eichler-Shimura Relations for Shimura Varieties of Hodge Type. ArXiv:2006.11745,

\bibitem{Levin-Weil-Restriction}B.~Levin, Local models for Weil-restricted groups. \emph{Comp. Math.}, \textbf{152} (2016), 2563--2601. 

\bibitem{Liu-HZ-Cycle}Y.~Liu, Hirzebruch-Zagier cycles and twisted triple product Selmer groups. \emph{Invent. Math.} \textbf{205} (2016) no.~3, 693--780.

\bibitem{LTXZZ-Rigidity}Y.~Liu and Y.~Tian and L.~Xiao and W.~Zhang and X.~Zhu, Deformation of rigid conjugate self-dual Galois representations. ArXiv:2108.06998.

\bibitem{LTXZZ-BBK}Y.~Liu and Y.~Tian and L.~Xiao and W.~Zhang and X.~Zhu, On the Beilinson-Bloch-Kato conjecture for Rankin-Selberg motives. \emph{Invent. Math} (2022).

\bibitem{Lourenco}J.~Louren\c{c}o, Grassmanniennes affines tordues sur les entiers. ArXiv:1912.11918.

\bibitem{Lusztig-q-analog}G.~Lusztig, Singularities, character formulas, and a $q$-analog of weight multiplicities. In \emph{Analysis and topology on singular spaces, II, III (Luminy, 1981)}, pp. 208--229, Ast\'erisque, 101-102, Soc. Math. France, Paris, 1983.


\bibitem{Lusztig-Unipotent}G.~Lusztig, Classification of unipotent representations of simple $p$-adic groups. \emph{Int. Math. Res. Not.} (1995), 517--589.

\bibitem{Mirkovic-Vilonen-GeoSat}I.~Mirkovi\'c and K.~Vilonen, Geometric Langlands duality and representations of algebraic groups over commutative rings. \emph{Ann. Math.} \textbf{166} (2007), 95--143.

\bibitem{Nekovar-Scholl-Plectic}J.~Nekov\'a\v{r} and A.J.~Scholl, Introduction to plectic cohomology. In \emph{Advances in the theory of automorphic forms and their L-functions}, Contemp. Math. Vol. 664, pp 321--337. Amer. Math. Soc., Providence, RI, 2016.

\bibitem{Ngo-Fundamental-Lemma}B.C.~Ng\^o, Le lemme fondamental pour les alg\`ebres de Lie. \emph{Publ. Math. l'IH\'ES} \textbf{111} (2010), 1--169.

\bibitem{Pappas-Rapoport-Ramified-Unitary}G.~Pappas and M.~Rapoport, Local models in the ramified case. III. Unitary groups.  \emph{J. Inst. Math. Jussieu} \textbf{8} (2009), 507--564.

\bibitem{PRS-Survey}G.~Pappas and M.~Rapoport and B.D.~Smithling, Local models of Shimura varieties, I. Geometry and combinatorics. In \emph{Handbook of Moduli, Vol. III}, pp. 135--217, International Press, Somerville, MA, 2013.

\bibitem{Pappas-Zhu-Local-Model}G.~Pappas and X.~Zhu, Local models of of Shimura varieties and a conjecture of Kottwitz.\emph{Invent. math.} \textbf{194} (2013), 147--254.

\bibitem{Scholze-Berkeley}P.~Scholze and J.~Weinstein, \emph{Berkeley lectures on $p$-adic geometry}. Ann. Math. Studies 207, Princeton University Press, 2020.

\bibitem{Tian-Xiao-Tate}
Y.~Tian and L.~Xiao, Tate cycles on quaternionic Shimura varieties over finite fields. \emph{Duke Math. J.} \textbf{168} (2019), 1551--1639.

\bibitem{Vollaard-Wedhorn-GU(n1)}I.~Vollaard and T.~Wedhorn, The supersingular locus of the Shimura variety of $GU(1,n-1)$, II. \emph{Invent. Math.} \textbf{184} (2011), 591--627.

\bibitem{Wiles-Fermat}A.~Wiles, Modular elliptic curves and Fermat's Last Theorem. \emph{Ann. Math.} \textbf{141} (1995), 443--551.

\bibitem{Wu-S=T}Z.~Wu, $S=T$ for Shimura Varieties and $p$-adic Shtukas. ArXiv:2110.10350.

\bibitem{Xiao-Zhu-Cycle}L.~Xiao and X.~Zhu, Cycles on Shimura varieties via geometric Satake. ArXiv:1707.05700.

\bibitem{Xiao-Zhu-Invariants}L.~Xiao and X.~Zhu, On vector-valued twisted conjugate invariant functions on a group. In \emph{Representations of Reductive Groups}, pp. 361--426, Proc. Symp. Pure Math. 101, Amer. Math. Soc., Providence, RI, 2019.

\bibitem{Xu-Zhu-Bessel}D.~Xu and X.~Zhu, Bessel $F$-isocrystals for reductive groups. In \emph{Invent. Math} (2022).

\bibitem{Xue-Finiteness}C.~Xue, Finiteness of cohomology groups of stacks of shtukas as modules over Hecke algebras, and applications. \emph{\'Epijournal de G\'eom\'etrie Alg\'ebrique} \textbf{4} (2020).

\bibitem{Xue-Smoothness}C.~Xue, Smoothness of cohomology sheaves of stacks of shtukas. ArXiv:2012.12833.

\bibitem{Yu-Integral-Satake}J.~Yu, The integral geometric Satake equivalence in mixed characteristic. ArXiv:1903.11132.

\bibitem{Yu-Jacquet-Langlands-Higher-Weights}J.~Yu, A geometric Jacquet-Langlands Transfer for automorphic forms of higher weights. In \emph{Trans. of Amer. Math. Soc.} (2022).

\bibitem{Yun-motive-exceptional-group}Z.~Yun, Motives with exceptional Galois groups and the inverse Galois problem. \emph{Invent. Math.} \textbf{196} (2014), 267--337. 

\bibitem{Yun-Zhu-Integral-Homology}Z.~Yun and X.~Zhu, Integral homology of loop groups via Langlands dual groups. \emph{Represent. Theory} \textbf{15}, (2011), 347--369.

\bibitem{Zhang-GGP}W.~Zhang, Fourier transform and the global Gan-Gross-Prasad conjecture for unitary groups. \emph{Ann. Math.} \textbf{180} (2014), 971--1049.

\bibitem{Zhu-S=T}X.~Zhu, $S=T$ for Shimura varieties. Preprint.

\bibitem{Zhu-Coherence-Conjecture}X.~Zhu, On the coherence conjecture of Pappas and Rapoport. \emph{Ann. Math.} \textbf{180} (2014), 1--85.

\bibitem{Zhu-Ramified-Sat}X.~Zhu, The Geometric Satake Correspondence for Ramified Groups. \emph{Ann. Sci. \'Ec. Norm. Sup\'er.}, \textbf{48} (2015),  409--451.

\bibitem{Zhu-CDM}X.~Zhu, Geometric Satake, categorical traces, and arithmetic of Shimura varieties. In \emph{Current Developments in Mathematics} (2016), 145--206.

\bibitem{Zhu-Satake-Mixed}X.~Zhu, Affine Grassmannians and the geometric Satake in mixed characteristic. \emph{Ann. Math.} \textbf{185} (2017), 403--492.

\bibitem{Zhu-FG-vs-HNY}X.~Zhu, Frenkel-Gross' irregular connection and Heinloth-Ng\^o-Yun's are the same. \emph{Selecta Math. (N.S.)} \textbf{23} (2017), no.~1, 245--274.

\bibitem{Zhu-PCMI}X.~Zhu, An introduction to affine Grassmannians and the geometric Satake equivalence. In \emph{IAS/Park City Math. Ser.} \textbf{24} (2017), 59--154.

\bibitem{Zhu-ICCM}X.~Zhu, A survey of local models of Shimura varieties. In \emph{Proceedings of the Sixth International Congress of Chinese Mathematicians}, pp. 247--265, Adv. Lectures in Math. Vol. 36, International Press , Boston, MA, 2017.

\bibitem{Zhu-Conjecture}X.~Zhu, Coherent sheaves on the stack of Langlands parameters. ArXiv:2008.02998.

\bibitem{Zhu-Integral-Satake}X.~Zhu, A note on integral Satake isomorphisms. ArXiv:2005.13056.



\end{thebibliography}
\end{document}